\documentclass[reqno,11pt]{amsart}
\pdfoutput=1

\usepackage[utf8]{inputenc}
\usepackage[italian,english]{babel}
\usepackage{amsfonts}
\usepackage{amssymb}
\usepackage{tikz}
\usepackage[right=3cm,bottom=2.5cm,footskip=30pt]{geometry}
\usepackage[autostyle]{csquotes}
\usepackage{color}
\usepackage{mathrsfs}
\usepackage{mathtools}
\usepackage{enumitem}
\usepackage{mathrsfs}
\setlist[enumerate,1]{leftmargin=.8cm}
\setlist[enumerate,2]{label={(\theenumi.\theenumii)},ref={(\theenumi.\theenumii)}}

\usepackage{subcaption}
\usepackage{tikz}
\usetikzlibrary{svg.path}
\definecolor{ColOrange}{HTML}{E27D60}
\definecolor{ColBlue}{HTML}{85CDCA}
\definecolor{ColYellow}{HTML}{E8A87C}
\definecolor{ColPink}{HTML}{C38D9D}
\definecolor{ColGreen}{HTML}{40B3A2}

\usepackage[toc]{appendix}
\usepackage{etoolbox}
\cslet{blx@noerroretextools}\empty
\usepackage[maxbibnames=9,giveninits=true,style=ext-alphabetic,sorting=nyt, backend=biber]{biblatex}
\DeclareFieldFormat{titlecase:title}{\MakeSentenceCase*{#1}}
\usepackage[colorlinks, citecolor=citegreen, linkcolor=refred,urlcolor=blue, unicode,psdextra,hypertexnames=false]{hyperref}
\usepackage{cleveref}
\crefname{enumi}{}{}
\crefname{enumii}{}{}
\expandafter\def\csname ver@etex.sty\endcsname{3000/12/31}

\usepackage{autonum}
\usepackage{gasymbols}

\newcommand{\tph}{\text{\tiny$(2)$}}

\definecolor{citegreen}{rgb}{0,0.3,0}
\definecolor{refred}{rgb}{0.5,0,0}

\usepackage[equation=section]{theorems}

\makeatletter
 \def\author@andify{
 \nxandlist {\unskip{} $\cdot$ \penalty-2}
 {\unskip {} $\cdot$ \penalty-2}
 {\unskip {} $\cdot$ \penalty-2}}
\makeatother

\newcommand{\orcid}[1]{\unskip {} \raisebox{-.3ex}{\href{https://orcid.org/#1}{\includegraphics{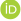}}}}

\definecolor{orcidlogocol}{HTML}{A6CE39}

\let\oldemail\email
\let\email\relax
\def\email#1{\oldemail{\href{mailto:#1}{\textcolor{black}{#1}}}}

\title[On the Isoperimetric Riemannian Penrose Inequality]{On the Isoperimetric Riemannian Penrose Inequality}

\author[L.~Benatti]{Luca Benatti\orcid{0000-0002-4685-7443}}
\address{L.~Benatti, Universit\`a degli Studi di Pisa,
Largo Bruno Pontecorvo 5, 56127 Pisa, Italy}
\email{luca.benatti@dm.unipi.it}

\author[M.~Fogagnolo]{Mattia Fogagnolo\orcid{0000-0002-5933-1344}}
\address{M.~Fogagnolo, Universit\`a di Padova, via Trieste 63, 35121 Padova (PD), Italy}
\email{mattia.fogagnolo@unipd.it}

\author[L.~Mazzieri]{Lorenzo Mazzieri\orcid{0000-0002-1597-8058}}
\address{L.~Mazzieri, Universit\`a degli Studi di Trento,
via Sommarive 14, 38123 Povo (TN), Italy}
\email{lorenzo.mazzieri@unitn.it}

\renewbibmacro*{url+urldate}{%
\iffieldundef{urlyear}
  {}
  {}
\usebibmacro{url}}

\renewcommand{\ncapa}{\mathfrak{c}}
\makeatletter
\newcommand{\svm}{s}%{\@ifnextchar_{g^s}{g_s}}
\makeatother
\newcommand{\eum}{\delta}
\newcommand{\ma}{\mathfrak{m}}

\begin{document}
 
\begin{abstract}

We prove that the Riemannian Penrose Inequality holds for Asymptotically Flat $3$-manifolds with nonnegative scalar curvature and connected horizon boundary, provided the optimal decay assumptions are met, which result in the $\ADM$ mass being a well-defined geometric invariant. Our proof builds on a novel interplay between the Hawking mass and a potential-theoretic version of it, recently introduced by Agostiniani, Oronzio and the third named author. 
As a consequence, we establish the equality between $\ADM$ mass and Huisken's Isoperimetric mass under the above sharp assumptions.
Moreover, we establish a Riemannian Penrose Inequality in terms of the Isoperimetric mass on any $3$-manifold with nonnegative scalar curvature, connected horizon boundary, and which supports a well-posed notion of weak Inverse Mean Curvature Flow. In particular, such Isoperimetric Riemannian Penrose Inequality does not require the asymptotic flatness of the manifold. The argument is based on a new asymptotic comparison result involving Huisken's Isoperimetric mass and the Hawking mass.

\end{abstract}
\maketitle


\medskip

\noindent \underline{\smash{Keywords}}: Penrose inequality, positive mass theorem, isoperimetric mass, inverse mean curvature flow, geometric inequalities.

\section{Introduction}
The concept of ADM mass, introduced in the early Sixties by the physicists Arnowitt, Deser and Misner \cite{arnowitt_coordinateinvarianceenergyexpressions_1961}, plays a central role in the modern development of mathematical relativity. This is due to its physical significance and its intriguing geometric properties. On the initial datum of a given isolated gravitational system, modelled by an Asymptotically Flat $3$-manifold $(M,g)$ with nonnegative scalar curvature, the $\ADM$ mass is given by
\begin{equation}\label{eq:mADM}
 \ma_{\ADM} = \lim_{j \to +\infty}\frac{1}{16\pi} \int_{\partial\Omega_j}g^{ij}(\partial_i g_{kj}- \partial_k g_{ij} ) \nu^k \dif \sigma ,
\end{equation}
where the flux integrals are taken along the boundaries of some exhaustion $(\Omega_j)_{j \in \N}$ of $M$, $\nu$ and $\dd \sigma$ being the unit normal vector field and the area measure on $\partial \Omega_j$ respectively.  In the above formula, the metric coefficients and their partial derivatives are expressed in a given \emph{asymptotically flat coordinate chart} $(x^1,x^2,x^3)$ (see \cref{def:asyflat}). Therefore, the value of the $\ADM$ mass is \textit{a priori} dependent on coordinates. However, in a couple of remarkable papers \cite{bartnik_massasymptoticallyflatmanifold_1986, chrusciel_boundaryconditionsspatialinfinity_1986}, it has been shown that the above expression defines a geometric invariant, provided the metric $g$ is $\CS^{1}_\tau$-Asymptotically Flat with $\tau > 1/2$. This assumption means that $\abs{g_{ij}-\delta_{ij}} = O_1(\abs{x}^{-\tau})$ at infinity. In this context, $\ma_{\ADM}$ is finite if and only if the scalar curvature is integrable on the whole space.

The most important statement about the $\ADM$ mass is the so-called Positive Mass Theorem, firstly proven by Schoen and Yau \cite{schoen_proofpositivemassconjecture_1979, schoen_proofpositivemasstheorem_1981} and subsequently by Witten \cite{witten_newproofpositiveenergy_1981,parker_wittenproofpositiveenergy_1982}. It says that $\ma_{\ADM}$ is nonnegative, and it is equal to zero if and only if $(M,g)$ is isometric to the flat Euclidean space. 

A sharper lower bound for the $\ADM$ mass was conjectured by Penrose \cite{penrose_nakedsingularities_1973} for Asymptotically Flat $3$-manifolds with nonnegative scalar curvature containing a horizon boundary $N$, which are known to model initial data sets of space-times containing a black hole. This inequality is expressed as
\begin{equation}
\sqrt{\frac{ \abs{ N}}{16 \pi}} \leq \ma_{\ADM},
\end{equation}
and the equality is achieved on the standard space-like slice of a Schwarzschild space-time. Loosely speaking, the total mass of an Asymptotically Flat solution containing a black hole must be at least as large as the total mass of a Schwarzschild solution with the same horizon's area. 

{Such a remarkable statement, known as the Riemannian Penrose Inequality, was first proven by Huisken and Ilmanen \cite{huisken_inversemeancurvatureflow_2001} for connected horizons, using the weak Inverse Mean Curvature Flow (IMCF for short). Almost at the same time, 
Bray~\cite{bray_proofriemannianpenroseinequality_2001} established a monotonicity formula for the ADM mass along a conformal flow of metrics, proving the same statement for multiple black holes. This proof was then extended by Bray and Lee~ \cite{bray_riemannianpenroseinequalitydimensions_2009} up to dimension 7. Recently, a new proof of the single black hole statement was proposed by the third author, in collaboration with Agostiniani, Mantegazza and Oronzio~\cite{agostiniani_riemannianpenroseinequalitynonlinear_2022}, using the level sets flow of $p$-harmonic functions.}

A common feature of the above-mentioned results is that the assumptions regarding the asymptotic behaviour of the metric are not optimal. Specifically, the decay rate required to establish the inequality is more stringent than what is needed to define the mass. More precisely, in Huisken-Ilmanen's proof through weak IMCF and similarly in the one through $p$-capacitary potentials \cite{agostiniani_riemannianpenroseinequalitynonlinear_2022}, the metric $g$ is assumed to be $\CS^1_1$-Asymptotically Flat. Moreover, the Ricci tensor is required to obey the following quadratic decay condition $\Ric \geq - \kst\abs{x}^{-2} g$, where $(x^1,x^2,x^3)$ are asymptotically flat coordinates and $\kst>0$ is a given constant. Indeed, the combination of these two conditions is sufficient to grant certain asymptotic properties of the functions whose level sets constitute the flow. These properties are then enough to compare the respective monotone quantities, i.e. the Hawking mass and its potential-theoretic counterpart, with the ADM mass. On the other hand, Bray's version of the Penrose Inequality \cite{bray_proofriemannianpenroseinequality_2001}, which also holds in the presence of possibly disconnected horizons, is proved requiring $\CS^{2}_\tau$-asymptotic flatness for some $\tau>1/2$. 

Our first main result confirms that the single horizon Riemannian Penrose Inequality holds whenever the sole decay conditions required to make the ADM mass a well-defined geometric invariant are met.

\medskip

\noindent
\emph{Throughout the whole article, except when explicitly stated, the Riemannian manifolds considered are agreed to be connected and with one single end.} 

\medskip

\begin{theorem}[Riemannian Penrose Inequality under optimal decay assumptions]
\label{thm:penrose_adm_intro}
Let $(M,g)$ be a complete $\CS^{1}_{\tau}$-Asymptotically Flat Riemannian $3$-manifold, $\tau >1/2$, with nonnegative scalar curvature and smooth, closed, outermost minimal boundary. Then, it holds
\begin{equation}\label{eq:adm_penrose}
 \sqrt{\frac{ \abs{N}}{16 \pi}} \leq \ma_{\ADM},
\end{equation}
where $N$ is any connected component of $\partial M$, with the equality satisfied if and only if $\partial M$ is connected and $(M,g)$ is isometric to a space-like slice of the Schwarzschild manifold
\begin{equation}
 \left(\R^3\smallsetminus \set{ \abs{x} < 2 \ma_{\ADM}},\, \left( 1+ \frac{ \ma_{\ADM}}{2 \abs{x}} \right)^{4}\eum_{ij}\,\dd x^i \otimes \dd x^j\right).
\end{equation}
\end{theorem}
A minimal boundary is called outermost if it is not shielded from infinity by some other minimal surface, in other words, if no other minimal surface homologous to $\partial M$ is present in $M\smallsetminus \partial M$.

The statement of \cref{thm:penrose_adm_intro} holds even when $\partial M = \varnothing$, where it provides the Positive Mass Theorem (PMT) under optimal decay assumption (see \cref{thm:pmt_adm} below). This version of the Positive Mass Theorem was previously known and implied by the one presented in \cite{lee_positivemasstheoremmanifolds_2015}.

\medskip

Along the route to \cref{thm:penrose_adm_intro} we establish a Riemannian Penrose Inequality for the \emph{Isoperimetric mass}, constituting our second main result. Huisken introduced this pretty general and very geometric concept of mass in \cite{huisken_isoperimetricconceptmassgeneral_2009}. We now recall its definition.

\begin{definition}[Isoperimetric mass] \label{def:isoperimetric_mass}
Let $(M,g)$ be a noncompact Riemannian $3$-manifold and let $\Omega$ be a bounded subset $\Omega\subset M$ with $\CS^{1,\alpha}$-boundary. The quasi-local Isoperimetric mass of $\Omega$ is defined as
\begin{equation}\label{eq:quasi-local_isoperimetric_mass_intro}
 \ma_{\iso}(\Omega)=\frac{2}{\abs{\partial \Omega}}\left( \abs{\Omega} -\frac{\abs{\partial \Omega}^{\frac{3}{2}}}{6 \sqrt{\pi}}\right).
\end{equation}
The Isoperimetric mass $\ma_{\iso}$ of $(M,g)$ is then defined as \begin{equation}\label{eq:isoperimetric_mass_intro}
 \ma_{\iso}=\sup_{(\Omega_j)_{j\in \N}}\limsup_{j \to +\infty} \ma_{\iso}(\Omega_j),
\end{equation}
where the supremum is taken among all exhaustions $(\Omega_j)_{j\in \N}$ consisting of domains with $\CS^{1,\alpha}$-boundary. 
\end{definition}
Huisken's definition in \cite{huisken_isoperimetricconceptmassgeneral_2009} involves families of sets with perimeters diverging at infinity rather than exhaustions. We remark here that all of our results hold for both notions of mass, and that they actually coincide once a global isoperimetric inequality is in force by \cite[Proposition 37]{jauregui_lowersemicontinuitymassconvergence_2019}.
It is immediately clear that the Isoperimetric mass is well-defined without needing any particular assumption on the asymptotic character of the metric. Hence, the question arises whether a Riemannian Penrose inequality is available for such a mass and for a broader class of initial data sets. An initial answer to this question can be found in \cite{jauregui_lowersemicontinuitymassconvergence_2019}. Following and strengthening the original argument of Huisken \cite{huisken_isoperimetricconceptmassgeneral_2009}, the authors prove that the Isoperimetric mass coincides with the $\ADM$ mass if the metric is assumed to be $\CS^2_\tau$-Asymptotically Flat at infinity. In the same setting, an alternative proof has been worked out in \cite[Appendix C]{chodosh_isoperimetryscalarcurvaturemass_2021}.

In the present paper, we extend the Isoperimetric Riemannian Penrose Inequality to the class of {\em strongly $1$-nonparabolic spaces}. 
This label is inspired by the concept of strongly nonparabolic spaces from \cite{ni_meanvaluetheoremsmanifolds_2007}. It grants that for every $\Omega\subset M$ with sufficiently regular boundary homologous to $\partial M$ there exists a proper solution $w_1\in \Lip_{\loc}(M\smallsetminus \Int \Omega)$ to the weak IMCF problem starting at $\Omega$
\begin{equation}\label{eq:IMCF}
 \begin{cases}
 \div\left( \dfrac{ \D w_1}{\abs{\D w_1}}\right) &=& \abs{ \D w_1} & \text{on $M\smallsetminus \overline{\Omega}$,}\\
 w_1 &=& 0 & \text{on $\overline{\Omega}$,}\\
 w_1 &\to& +\infty & \text{as $\dist(x,o) \to +\infty$}.
 \end{cases}
\end{equation} 
We can now state our second main result.

\begin{theorem}[Isoperimetric Riemannian Penrose Inequality]
\label{thm:isoperimetric_penrose_intro}
Let $(M,g)$ be a complete strongly $1$-nonparabolic Riemannian $3$-manifold, with nonnegative scalar curvature and smooth, compact, outermost minimal boundary. Then, it holds
\begin{equation}\label{eq:isoperimetric_penrose}
 \sqrt{\frac{ \abs{N}}{16 \pi}} \leq \ma_{\iso},
\end{equation}
where $N$ is any connected component of $\partial M$, with the equality satisfied if and only if $\partial M$ is connected and $(M,g)$ is isometric to a space-like slice of the Schwarzschild manifold
\begin{equation}
 \left(\R^3\smallsetminus \set{ \abs{x} < 2 \ma_{\iso}} ,\, \left(  1+ \frac{ \ma_{\iso}}{2 \abs{x}} \right)^{4}\eum_{ij}\, \dd x^i \otimes \dd x^j\right).
\end{equation}
\end{theorem}
 We notice once more that the above result does not require any specific decay of the metric coefficients towards the Euclidean model. The {\em strong $1$-nonparabolicity condition} is satisfied for example on $3$-manifolds that merely support a Euclidean-like Isoperimetric inequality \cite{xu_isoperimetrypropernessweakinverse_2023}. If the boundary is not present or not outermost, one can still draw out of the arguments leading to \cref{thm:isoperimetric_penrose_intro} an Isoperimetric Positive Mass Theorem, see \cref{thm:isoperimetric_pmt}.

In deriving \cref{thm:penrose_adm_intro} from \cref{thm:isoperimetric_penrose_intro}, we also settle the following result about the equivalence of masses.

\begin{theorem}\label{thm:mass_equivalence}
Let $(M,g)$ be a complete $\CS^{1}_{\tau}$-Asymptotically Flat Riemannian $3$-manifold, $\tau >1/2$, with nonnegative scalar curvature and possibly empty smooth, compact and minimal boundary. Then, 
\begin{equation}\label{eq:equivalence_mass}
\ma_{\ADM}=\ma_{\iso}.
\end{equation}
\end{theorem}

This result should be viewed as an improvement on \cite[Theorem 3]{jauregui_lowersemicontinuitymassconvergence_2019}.
 As noted by Chruściel, \cref{eq:equivalence_mass} provides an alternative proof of the geometric invariance of the ADM mass in $3$-manifolds of nonnegative scalar curvature and in the optimal decay regime. We mention that our results never require that the masses are finite. In particular, \cref{eq:equivalence_mass} also infers that the $\ADM$ mass is finite if and only if the Isoperimetric mass is finite.

\medskip

All the above-mentioned results rely on a new asymptotic comparison argument, which compares the Hawking mass along the weak IMCF $w_1$ with the quasi-local Isoperimetric mass of the evolving regions. The crux of this analysis is  \cref{thm:mass_control_IMCF}.
Using this result and bringing to light a new interplay between the Hawking mass $\ma_{H}$ and its potential theoretic version $\ma^{\tph}_{H}$ (see \cref{eq:2hawking}), we pave the way to the proof of \cref{thm:penrose_adm_intro}. In this context, the quantity $\ma^{\tph}_{H}$ was recently introduced in \cite{agostiniani_greenfunctionproofpositive_2024}. The authors also show its monotonicity along the level set flow of the harmonic capacitary potential of the initial set $\partial\Omega$, mirroring the monotonic behaviour of the Hawking mass along the weak IMCF. 

To obtain the desired asymptotic comparison between the Hawking mass and the $\ADM$ mass, Huisken-Ilmanen \cite{huisken_inversemeancurvatureflow_2001} had to require extra assumptions on the behaviour of the metric. These are needed to compensate for the lack of sufficiently refined asymptotics for the weak IMCF. Ultimately building on Green's representation formula, the potential theoretic approach rather allows to derive very precise second-order asymptotic expansions for solutions of the Laplace equation, which leads to control of the mass $\ma^{\tph}_{H}$ in terms of the ADM mass. This approach permits giving a finite bound on the Hawking mass and also provides some new insights about the asymptotic behaviour of the weak IMCF. Building on it, one can finally provide the sharp upper bound for the Hawking mass, which has two main consequences. First, it directly yields \cref{thm:penrose_adm_intro} using $\partial \Omega =N$ as the starting set for the weak IMCF. Second, exploiting Jauregui-Lee's techniques \cite{jauregui_lowersemicontinuitymassconvergence_2019} (see \cref{thm:bound_hawking_to_bound_iso} below) we infer $\ma_{\iso} \leq \ma_{\ADM}$. This inequality is enough to ensure \cref{thm:mass_equivalence} since the reverse one is already contained in  \cite{fan_largespheresmallspherelimitsbrownyork_2009}. 

\medskip

We finally observe that
\cref{thm:isoperimetric_penrose_intro} provides a new piece of information also in the basic situation of manifolds with nonnegative scalar curvature that are $\CS^0$-asymptotic to a Schwarzschild space of nonnegative mass $\ma$ (see \cref{def:asysvar}). Indeed, in \cref{thm:asysvarchild} below, we show that such parameter satisfies the expected Penrose Inequality and coincides in fact with the Isoperimetric mass. We emphasise that in this setting the classical notion of $\ADM$ mass is not even clearly defined, as it would involve the quantitative decay of the derivatives of the metric coefficients codified in the $\CS^1_\tau $-asymptotic condition, with $\tau > 1/2$. 

\begin{theorem}[Penrose Inequality in Asymptotically Schwarzschildian manifolds]
\label{thm:asysvarchild}
Let $(M, g)$ be a  $\CS^0_\tau$-Asymptotically Schwarzschildian of mass $\ma \geq 0$, $\tau>0$, with nonnegative scalar curvature and possibly empty smooth, closed and outermost minimal boundary.  Then, it holds
 \begin{equation}
 \label{eq:penrose-svarchild}
 \sqrt{\frac{\abs{N}}{16 \pi}}\leq \ma,
\end{equation}
where $N$ is any connected component of $\partial M$, with the equality satisfied if and only if $\partial M$ is connected and $(M,g)$ is isometric to a space-like slice of the Schwarzschild manifold
\begin{equation}
 \left(\R^3\smallsetminus \set{ \abs{x} < 2 \ma} , \, \left( 1+ \frac{ \ma}{2 \abs{x}} \right)^{4}\eum_{ij}\, \dd x^i  \otimes \dd x^j\right).
\end{equation}
Moreover, we have $\ma = \ma_{\iso}$. 
\end{theorem}

\subsection*{Perspectives} 

Back in \cite{agostiniani_riemannianpenroseinequalitynonlinear_2022}, a family of Hawking-type masses defined along the level sets of $p$-capacitary potentials, $p \in (1, 3)$, was introduced and proved to be monotone. In \cite{benatti_nonlinearisocapacitaryconceptsmass_2023a}, we use such $p$-Hawking masses and the monotone quantities recently obtained in \cite{chan_monotonicitygreenfunctions_2022} to prove new Positive Mass Theorems and Penrose-type inequalities for Iso-$p$-capacitary masses. These masses interpolate between Huisken's Isoperimetric mass and the capacitary mass, proposed by Jauregui \cite{jauregui_admmasscapacityvolumedeficit_2020}. 

\smallskip

We believe that the techniques presented here, that allow us to work with very mild asymptotic assumptions, could have implications for a conjecture of Huisken \cite[p. 2221-2223]{cederbaum_mathematicalaspectsgeneralrelativity_2021}. This conjecture states that a Positive Mass Theorem for the Isoperimetric mass should hold on Riemannian manifolds with $\CS^0$-metrics and a suitable weak definition of nonnegative scalar curvature. Versions of this conjecture have been proved in \cite{antonelli_positivemass_2024, jauregui2024notehuiskensisoperimetricmass} for Asymptotically Flat manifolds. Another consequence of such papers is that the equivalence between masses stated in \cref{thm:asysvarchild} fails in Schwartzschild $3$-manifolds with negative mass.
Remaining in the realm of spaces with low regularity, Burkhardt-Guim \cite{burkhardt-guim_admmassmetricsdistortion_2022} settled the basic properties of a new $\CS^0$-concept of $\ADM$ mass when a Ricci-flow-related notion of nonnegative scalar curvature is assumed. It would be interesting to investigate the relations between Burkhardt-Guim's notion and the Isoperimetric mass.

\subsection*{Outline of the paper} 
The paper is structured as follows. In \cref{sec:isoperimetric}, we first recall the definition of the weak Inverse Mean Curvature Flow and its main features in an ambient $3$-manifold with nonnegative scalar curvature. Then, we prove the Asymptotic Comparison Lemma \ref{thm:mass_control_IMCF}, which constitutes the main new ingredient in the proof of the Isoperimetric Riemannian Penrose Inequality. 
We discuss some relevant topological properties of  $3$-manifolds with nonnegative scalar curvature from the perspective of 
mean-convex mean-curvature flow with surgeries  \cite{brendle_meancurvatureflowsurgery_2018}. We proceed with the 
proof of  \cref{thm:isoperimetric_penrose_intro}.
In \cref{sec:asysvar}, we provide an existence result for the weak IMCF in $\CS^0$-Asymptotically Flat Riemannian manifolds and describe its asymptotic features. Then, we prove \cref{thm:asysvarchild}. \cref{sec:admpenrose} is devoted to the proof of \cref{thm:penrose_adm_intro,thm:mass_equivalence}. We first work out the asymptotic analysis for the capacitary potential to get a bound from above for the $2$-Hawking mass \cref{eq:2hawking}. Then, we apply it jointly with new fine asymptotic estimates along the weak IMCF to sharply control the Hawking mass and complete the proof. The manuscript ends with two Appendices. In the first one, a general version of the de l'H\^opital's rule is stated and proved, as it constitutes a fundamental ingredient in the proof of \cref{thm:isoperimetric_penrose_intro}. In the second one, we provide some estimates in nonlinear potential theory that we need to study the IMCF in $\CS^0$-Asymptotically Flat Riemannian manifolds.

\subsection*{Acknowledgements}
A substantial part of this work has been carried out during the authors' attendance to the \emph{Thematic Program on Nonsmooth Riemannian and Lorentzian Geometry} that took place at the Fields Institute in Toronto in 2022. The authors warmly thank the staff, the organizers and the colleagues for the wonderful atmosphere and the excellent working conditions set up there. L. B. is supported by the European Research Council’s (ERC) project n.853404 ERC VaReg -- \textit{Variational approach to the regularity of the free boundaries}, financed by the program Horizon 2020, by PRA\_2022\_11 and by PRA\_2022\_14.
M. F. is supported by the Project funded under the National Recovery and Resilience Plan (NRRP), Mission 4 Component 2 Investment 1.5 - Call for tender No. 3277 of 30 dicembre 2021 of Italian Ministry of University and Research funded by the European Union – NextGenerationEU; Project code: ECS00000043, Concession Decree No. 1058 of June 23, 2022, CUP C43C22000340006, Project title ``iNEST: Interconnected Nord-Est Innovation Ecosystem''.
During the preparation of the work M. F. was funded by the European Union – NextGenerationEU and by the University of Padova under the 2021 STARS Grants@Unipd programme ``QuASAR". The authors are members of Gruppo Nazionale per l’Analisi
Matematica, la Probabilit\`a e le loro Applicazioni (GNAMPA), which is part of the Istituto
Nazionale di Alta Matematica (INdAM), and are partially funded by the GNAMPA project ``Problemi al bordo e applicazioni geometriche".  The authors are grateful to V. Agostiniani, G. Antonelli, S. Borghini, A. Carlotto, A. Malchiodi, F. Oronzio, A. Pluda and M. Pozzetta for their interest in the work and for pleasureful and useful conversations on the subject. The authors thank the anonymous referee for their important comments and remark, that concurred to improve the quality of the paper.

\section{Isoperimetric Riemannian Penrose Inequality}
\label{sec:isoperimetric}
This section aims to prove the Isoperimetric Penrose Inequality \cref{thm:isoperimetric_penrose_intro}. The main ingredients are two. The first one is the Geroch monotonicity formula.
Proved by Huisken-Ilmanen in \cite[Geroch Monotonicity Formula 5.8]{huisken_inversemeancurvatureflow_2001}, it states that $t\mapsto \ma_{H}(\partial \Omega_t)$ is monotone nondecreasing, where 
\begin{equation}\label{eq:hawking-mass}
 \ma_{H}(\partial \Omega)= \frac{ \abs{ \partial \Omega}^{\frac{1}{2}}}{16 \pi^{\frac{3}{2}}} \left( 4 \pi - \int_{\partial \Omega}\frac{\H^2}{4} \dif \sigma \right)
\end{equation}
is the \textit{Hawking mass} \cite{hawking_gravitationalradiationexpandinguniverse_1968} and $\Omega_t =\set{w_1 \leq t}$ are the sublevels of the weak IMCF \cref{eq:IMCF}. The second one is a novel, fairly simple, asymptotic comparison between $\ma_H(\partial\Omega_t)$ and $\ma_{\mathrm{iso}}(\Omega_t)$ as $t\to +\infty$ (see \cref{thm:mass_control_IMCF}). Indeed, we can apply the de l'H\^opital's rule to differentiate the quasi-local Isoperimetric mass. Exploiting some favourable properties of the weak IMCF, we manipulate the resulting expression to recover the Hawking mass.

\subsection{The weak Inverse Mean Curvature Flow}
The theory of weak IMCF has been introduced in \cite{huisken_inversemeancurvatureflow_2001} as the main tool for proving the Riemannian Penrose Inequality. In this work, we will not focus on the existence theory of the weak IMCF. We consider Riemannian manifolds $(M,g)$ where any bounded $\Omega$ with a sufficiently smooth boundary homologous to $\partial M$ admits a weak proper solution to \cref{eq:IMCF}. 
In the absence of a boundary, where the Riemannian Penrose inequality consists in the (Riemannian) Positive Mass Theorem, we are going to exploit IMCF issuing from a point, which is
\begin{equation}\label{eq:IMCFfromapoint}
 \begin{cases}
 \div\left( \dfrac{ \D w_1}{\abs{\D w_1}}\right) &=& \abs{ \D w_1} & \text{on $M\smallsetminus \set{o}$,}\\
 w_1(x) &\to& -\infty & \text{as $\dist(x,o) \to 0$,}\\
 w_1(x) &\to& +\infty & \text{as $\dist(x,o) \to +\infty$}.
 \end{cases}
\end{equation}
We refer to \cite{huisken_inversemeancurvatureflow_2001} for the precise definition of weak solution to the IMCF.

\medskip

We give our definition of strong $1$-nonparabolicity. 

\begin{definition}[strong $1$-nonparabolicity]
\label{def:strongly1nonpar}
Let  $(M,g)$ be a Riemannian manifold with possibly empty compact boundary $\partial M$. 
\begin{enumerate}
\item If $\partial M \neq \varnothing$, we say that $(M, g)$ is strongly $1$-nonparabolic if there exists a proper weak solution $w_1$ to \cref{eq:IMCF} starting at $ \partial M$.
\item If $\partial M = \varnothing$, we say that $(M, g)$ is strongly $1$-nonparabolic if, for any $o \in M$, there exists a proper weak solution $w_1$ to \cref{eq:IMCFfromapoint}.
\end{enumerate}
\end{definition}

The above is in accordance with Ni's definition of \emph{strongly nonparabolic} \cite{ni_meanvaluetheoremsmanifolds_2007}. We recall that a Riemannian manifold is strongly nonparabolic if it admits a positive Green's function converging to zero at infinity. 
The prefix $1$ is motivated by the fact the IMCF can be constructed as suitable limit as $p\to 1^+$ of $p$-harmonic functions \cite{moser_inversemeancurvatureflow_2007,moser_inversemeancurvatureflow_2008,kotschwar_localgradientestimatesharmonic_2009,mari_flowlaplaceapproximationnew_2022}.

\smallskip

Without mentioning it, we use that a strongly $1$-nonparabolic manifold admits a solution to \cref{eq:IMCF} starting at any $\Omega$ with sufficiently regular boundary homologous to $\partial M$. This is granted by \cite[Weak Existence Theorem 3.1]{huisken_inversemeancurvatureflow_2001}.

\begin{remark}
By \cite[Theorem 1.2]{xu_isoperimetrypropernessweakinverse_2023}, every $(M,g)$ of dimension $3 \leq n \leq 7$ admitting a Euclidean-like Isoperimetric Inequality is strongly $1$-nonparabolic. Actually, the author provides the existence of the weak IMCF only starting at a given hypersurface. However, building on the quantitative diameter estimate for minimising hulls proved there, one can also deduce the existence of a weak IMCF issuing from any given point.
\end{remark}

In the following lemma, we recall some geometric properties that hold along the level sets of the weak IMCF (see \cite[Regularity Theorem 1.3, Minimizing Hull Property 1.4 and Exponential Growth Lemma 1.6]{huisken_inversemeancurvatureflow_2001}).

\begin{lemma}
Let $(M,g)$ be a Riemannian $3$-manifold possibly with boundary. Let $w_1$ be the solution to \cref{eq:IMCF} starting at some bounded subset $\Omega \subseteq M$. 
\begin{enumerate}
    \item $\Omega_t = \set{w_1\leq t}$ is strictly outward minimising with $\CS^{1,\alpha}$-boundary,
    \item $\abs{\partial \Omega_t}= \ee^{t}\abs{\partial \Omega^*}$.
\end{enumerate}
\end{lemma}
 We observe that the map $t \mapsto \abs{\set{w_1 \leq t }}$ is monotone increasing, hence it has countably many discontinuities. Out of these jump times, one has $\partial \set{w_1 \leq t } = \partial \set{w_1<t} = \set{ w_1 =t} $. We also recall that $\partial \Omega_t$ admits a weak notion of mean curvature $\H$ and that
 \begin{align}
    \H = \abs{\nabla w_1}>0 && \text{a.e. on } \partial \Omega_t \qquad \text{ for a.e. } t,
 \end{align}
 as stated in \cite[(1.12) and Lemma 5.1]{huisken_inversemeancurvatureflow_2001}.

The Hawking mass \cref{eq:hawking-mass} is well-defined at almost every level of the weak IMCF, as proved in \cite[Section 1, pp. 369-371]{huisken_inversemeancurvatureflow_2001}. A fundamental property of the Hawking mass consists on its monotonicity along the weak IMCF \cite[Geroch Monotonicity Formula 5.8]{huisken_inversemeancurvatureflow_2001}. We report the statement of this theorem for ease of future reference. In the following, as well as in the rest of the paper, we denote with $\h$ and $\mathring{\h}$ the weak second fundamental form and its traceless part.
\begin{theorem}[Monotonicity of the Hawking mass]\label{thm:Geroch_monotonicity_formula}
    Let $(M,g)$ be a complete, noncompact Riemannian $3$-manifold with nonnegative scalar curvature and possibly with a smooth and closed boundary. Assume that $H_2(M,\partial M; \Z)=\set{0}$. Let $\Omega \subseteq  M$ be with connected $\CS^1$-boundary homologous to $\partial M$ such that $\h \in L^2(\partial \Omega)$. Let $w_1 \in \Lip_{\loc}(M \smallsetminus \Int \Omega)$ be a solution to the problem \cref{eq:IMCF}. Then, the function $t\mapsto \ma_{H}(\partial \Omega_t)$ is a monotone nondecreasing $\BV_{\loc}(0,+\infty)$ function and
    \begin{equation}\label{eq:geroch-monotonicity}
        \frac{\dd}{\dd t}\ma_{H}(\partial \Omega_t) \geq\frac{\sqrt{ \abs{ \partial \Omega_t}}}{(4 \sqrt{\pi})^3} \int_{\partial \Omega_t}2 \frac{\abs{\D^{\top} \H}^2}{\H^2} + \vert{ \mathring{\h}}\vert^2 \dif \sigma
    \end{equation}
    holds in the weak sense, where $\Omega_t= \set{w_1 \leq t }$.
\end{theorem}
The condition on the second homology group is added to ensure that the boundaries of the level sets of $w_1$ remain connected, which is essential for establishing \cref{eq:geroch-monotonicity}. Indeed, the proof of this inequality exploits the Gauss-Bonnet Theorem to control $8\pi - \int_{\partial \Omega_t} \mathrm{R}^\top \dif \sigma$, where $\mathrm{R}^\top$ is the scalar curvature of the induced metric on $\partial \Omega_t$. 
This assumption is always satisfied in the setting of \cref{thm:penrose_adm_intro,thm:mass_equivalence,thm:asysvarchild} (see \cref{thm:topology} below). The relation between the second homology group and the boundaries of the level sets of $w_1$ is tackled in the following proposition.
For the sake of completeness, we provide the proof by Bray-Miao \cite[Proposition 1]{bray_capacitysurfacesmanifoldsnonnegative_2008}, including a few more details.

\begin{proposition}\label{prop:connectedness_of_boundaries}
Let $(M, g)$ be a strongly $1$-nonparabolic Riemannian $3$-manifold satisfying $H_2(M, \partial M; \Z) =\set{0}$. Then, if $w_1$ is a weak IMCF \cref{eq:IMCF} starting at $\Omega$ with connected boundary homologous to $\partial M$, then the surfaces $\partial\set{w_1 <t}$ and $\partial\set{w_1 \leq t}$ are connected for any $t > 0$.
\end{proposition}
\begin{proof}
It suffices to prove the claim out of the countably many $t$'s at which the flow jumps. In particular, $\partial\set{w_1 <t} = \partial\set{w_1 \leq t} = \set{w_1 = t}$ is satisfied at these $t$'s. Thus, by \cite[(1.10)]{huisken_inversemeancurvatureflow_2001}, we can approximate $\partial\set{w_1 <t}$
and $\partial\set{w_1 \leq t}$ in $\CS^1$ by sets $\set{w_1 = s_j}$ for sequences $s_j \to t^-$ and $s_j \to t^+$ respectively. 

\smallskip

Let $\Sigma_t$ be a connected component of $\set{w_1 = t}$. By the assumption $H_2(M,\partial M; \Z)=\set{0}$, there exists a bounded open set $D_t$ such that $\partial D_t = \Sigma_t$. Here, $\partial D_t$ denotes the topological boundary of $D_t$ in the manifold $M$. We are committed to show that $D_t = \set{w_1 < t}$. This assertion implies that $\set{w_1=t}$ is connected, as it would coincide with $\partial D_t =\Sigma_t$. While a \textit{priori} $D_t$ may or may not contain $\partial M$, the result also implies that $\partial M \subseteq D_t$.

We first show that $D_t\subseteq \set{w_1 \leq t}$.
By contradiction, suppose there were $x \in D_t$ with $w_1(x) > t$. Take $y \in M \smallsetminus \overline{D}_t$ such that $w_1(y) > t$. This point $y$ exists because $\overline{D_t}$ is compact and $w_1 \to + \infty$ at infinity. As shown in the proof of \cite[Connectedness Lemma 4.2 (ii)]{huisken_inversemeancurvatureflow_2001}, $\set{w_1 > t}$ is connected since $(M,g)$ has only one end. Then, we can join $x$ and $y$ with a continuous path, entirely lying in $\set{w_1 > t}$. This curve must cross $\partial D_t \subset \set{w_1 = t}$, which is a contradiction. 

We now show that ${D}_t \subseteq \set{w_1 < t}$.
By contradiction, assume that there is $x \in D_t$ satisfying $w_1(x) = t$.
As recalled above, we can approximate $\Sigma_t$ in $\CS^1$ by a connected component $\Sigma_s$ of $\set{w_1 = s}$ as $s \to t^-$. Then, for $s < t$ sufficiently close to $t$, the point $x$ would belong to $D_s$ but not to $\set{w_1 \leq s}$, yielding a contradiction to what just proved.

Finally, we show that $\set{w_1 < t} \subseteq {D_t}$, completing the proof. By contradiction, assume that there is $x \in M \smallsetminus {D_t}$ at which $w_1 (x) < t$. 
Then, $x \in M \smallsetminus \overline{D}_t$ since $w_1(x) = t$ on $\partial D_t$. Consider now a point $y \in D_t$. We proved that $w_1(y) < t$. The proof of \cite[Connectedness Lemma 4.2 (ii)]{huisken_inversemeancurvatureflow_2001} shows that if $\partial \Omega$ is connected then $\set{w_1 < t}$ is connected as well. Again, we can join $x$ and $y$ with a continuous path in $\set{w_1 <t}$. Since the curve must cross $\partial D \subset \set{w_1 = t}$, we have the desired contradiction. 
\end{proof}

\smallskip

Take an initial set $\Omega\subseteq M$ with connected boundary homologous to $\partial \Omega$. If one can produce an asymptotic bound for the Hawking mass, the monotonicity transfers this same bound to the Hawking mass of $\partial \Omega$. Observe that a single component $N$ is not homologous to $\partial M$ when the latter consists of multiple connected components. Thus, we cannot rely on \cref{thm:Geroch_monotonicity_formula} to get a bound on $\ma_H(N)$ from an asymptotic estimate. However, \cite[Geroch Monotonicity (Multiple Boundary Components) 6.1]{huisken_inversemeancurvatureflow_2001} implies that having an upper bound for subsets with boundary homologous to $\partial M$ is actually enough to control the Hawking mass of any component of the horizon. We summarize the procedure in the following statement.

\begin{proposition}\label{prop:jumping_to_homologous}
    Let $(M,g)$ be a complete, strongly $1$-nonparabolic Riemannian $3$-manifold with nonnegative scalar curvature and with smooth and closed boundary. Assume that the $\partial M$ consists of embedded minimal spheres and $H_2(M,\partial M; \Z)=\set{0}$. Let $E \subseteq  M$ be a bounded subset with connected smooth boundary, such that any connected component of $\partial M$ is either contained in $E$ or disjoint from $E$. Then, there exists $\Omega \subseteq  M$ with the following properties.
    \begin{enumerate}
        \item $\Omega$ is strictly outward minimising with connected $\CS^{1,\alpha}$-boundary homologous to $\partial M$ and with $\h\in L^2(\partial \Omega)$ (see \cite[Section 5, pp. 401-402]{huisken_inversemeancurvatureflow_2001}).
        \item $\ma_{H}(\partial E) \leq \ma_{H}( \partial \Omega)$.
        \item There exists a solution $w_1\in \Lip_{\loc}(M\smallsetminus \Int \Omega) $  to \cref{eq:IMCF} starting at $\Omega$. 
    \end{enumerate}
\end{proposition}

The condition on the topology of the boundary ensures that the manifold $M'$ obtained by filling in the boundary with balls satisfies $H_2(M';\Z)= \set{0}$ and $\partial E$ is homologous to $\partial M' = \varnothing$. Hence, \cref{prop:connectedness_of_boundaries} applies. The level sets of the weak IMCF remain connected before touching a new connected component of $\partial M$. Then one can jump over the horizon and restart the flow at a new connected boundary with greater Hawking mass, as described in \cite[Geroch Monotonicity (Multiple Boundary Components) 6.1]{huisken_inversemeancurvatureflow_2001}.

\smallskip

Building on the original argument of Huisken \cite{huisken_isoperimetricconceptmassgeneral_2009}, Jauregui and Lee \cite{jauregui_lowersemicontinuitymassconvergence_2019} showed that the isoperimetric mass and the $\ma_{\ADM}$ do coincide under the assumptions of \cite{huisken_inversemeancurvatureflow_2001}.
The crucial step in their proof consists in showing that, if there exists $\zeta$ such that 
\begin{equation}
\label{eq:boundhawking-discussion}
\ma_{H} (\partial \Omega) \leq \zeta
\end{equation}
for any $\Omega$ belonging to a suitable class, then one can in fact derive $\ma_{\iso} \leq \zeta$.
Through a perusal of the proof of \cite[Theorem 17]{jauregui_lowersemicontinuitymassconvergence_2019}, it is easily checked that this step only involves the $\CS^0$-character of the metric. 
The authors assume $\CS^2_1$-asymptotic flatness to apply \cite{huisken_inversemeancurvatureflow_2001}, allowing \cref{eq:boundhawking-discussion} to hold with $\zeta = \ma_{\ADM}$. However, our purposes require the full strength of Jauregui-Lee's analysis, which we summarise in the following statement.

\begin{theorem}\label{thm:bound_hawking_to_bound_iso}
        Let $(M,g)$ be a complete $\CS^0$-Asymptotically Flat Riemannian $3$-manifold with a possibly empty smooth minimal boundary. Assume $\ma_{H} (\partial \Omega) \leq \zeta$ for any bounded subset with connected smooth boundary and such that any connected component of $\partial M$ is either contained in $\Omega$ or disjoint from $\Omega$. Then,  $\ma_{\iso} \leq \zeta$.
\end{theorem}

\subsection{An asymptotic comparison Lemma.}
The following lemma contains the most relevant asymptotic estimate of the present work. Used in combination with Huisken-Ilmanen's proof of the Geroch Monotonicity formula \cref{thm:Geroch_monotonicity_formula}, it leads directly to the Isoperimetric Riemannian Penrose Inequality \cref{thm:isoperimetric_penrose_intro}.

\begin{lemma}[Asymptotic Comparison Lemma]
\label{thm:mass_control_IMCF}
Let $(M,g)$ be a complete, noncompact Riemannian $3$-manifold with nonnegative scalar curvature and possibly with smooth and closed boundary. Assume that {$H_2(M,\partial M; \Z)=\set{0}$}. Let $\Omega \subseteq  M$ be bounded with connected $\CS^1$-boundary homologous to $\partial M$ satisfying $\h\in L^2(\partial \Omega)$ and admitting a solution $w_1 : M \smallsetminus \Int \Omega\to \R$ to the problem \cref{eq:IMCF}. Then
\begin{equation}\label{eq:main_inequality_IMCF}
    \lim_{t \to +\infty} \ma_H(\partial \Omega_t)\leq \liminf_{t \to +\infty} \ma_{\iso}(\Omega_t),
\end{equation}
where $\Omega_t = \set{ w_1 \leq t}$. 
\end{lemma}

\begin{proof} Suppose that the right-hand side in \cref{eq:main_inequality_IMCF} is finite otherwise there is nothing to prove. By the coarea formula, we have
\begin{equation}
    \abs{\Omega_t \smallsetminus \Crit w_1} = \int_0^t \int_{\partial \Omega_s} \frac{1}{\H} \dif \sigma \dif s,
\end{equation}
which ensures that $t\mapsto \abs{\Omega_t\smallsetminus \Crit w_1}$ is locally absolutely continuous and monotone and that $t\mapsto \int_{\partial \Omega_t} 1/\H \dif \sigma$ is a well-defined $L^1_{\loc}$ function. We are therefore in position to apply the de L'H\^opital's rule (\cref{thm:delhopital}) to $t\mapsto \ma_{\iso}( \Omega_t)$, obtaining
\begin{align}\label{eq:IMCF_de_lhopital}
    \begin{split}
        \liminf_{t \to +\infty} \ma_{\iso}(\Omega_t)&\geq \liminf_{t \to +\infty} \frac{2}{\abs{\partial \Omega_t}}\left( \abs{\Omega_t\smallsetminus\Crit w_1} -\frac{\abs{\partial \Omega_t}^{\frac{3}{2}}}{6 \sqrt{\pi}}\right)\\&\geq\liminf_{t \to +\infty} \frac{2}{\abs{\partial \Omega_t}}\left(\kern.1cm \int_{\partial \Omega_t} \frac{1}{\H} \dif \sigma -\frac{\abs{\partial \Omega_t}^{\frac{3}{2}}}{4 \sqrt{\pi}}\right).
    \end{split}
\end{align}
The H\"{o}lder's inequality yields
\begin{equation}
    \int\limits_{\partial \Omega_t }\frac{1}{\H} \dif \sigma \geq \abs{\partial \Omega_t}^{\frac{3}{2}}\left( \kern.1cm \int\limits_{\partial \Omega_t} \H^2\dif \sigma \right)^{-\frac{1}{2}},
\end{equation}
which gives
\begin{equation}\label{eq:second_IMCF_de_lhopital}
\begin{split}
    \liminf_{t \to +\infty}\ma_{\iso}(\Omega_t) &\geq \liminf_{t \to +\infty}2 \left(\frac{\abs{ \partial \Omega_t}}{\int_{\partial \Omega_t}\H^2 \dif \sigma}\right)^{\frac{1}{2}}\left( 1- \frac{1}{4\sqrt{\pi}}\left(\kern.1cm\int_{\partial \Omega_t} \H^2 \dif \sigma\right)^{\frac{1}{2}}\right)\\
    &=\liminf_{t \to +\infty}2 \left(\frac{\abs{ \partial \Omega_t}}{\int_{\partial \Omega_t}\H^2 \dif \sigma}\right)^{\frac{1}{2}}\frac{ 1- \frac{1}{16\pi}\int_{\partial \Omega_t} \H^2 \dif \sigma}{1+( \frac{1}{16\pi}\int_{\partial \Omega_t} \H^2 \dif \sigma)^{\frac{1}{2}} }\\
    &=\liminf_{t \to +\infty}\frac{2 \ma_H(\partial \Omega_t)}{( \frac{1}{16\pi}\int_{\partial \Omega_t} \H^2 \dif \sigma)^{\frac{1}{2}}+\frac{1}{16\pi}\int_{\partial \Omega_t} \H^2 \dif \sigma} \, .
\end{split}
\end{equation}
The statement follows once the following claim is proved.
\begin{claim}
There exists a divergent increasing sequence $(t_n)_{n \in \N}$ realising the limit inferior in the right-hand side of \cref{eq:second_IMCF_de_lhopital} and such that $\int_{\partial \Omega_{t_n}} \H^2 \dif \sigma\to 16 \pi$ as $n\to +\infty$.
\end{claim}
\noindent Indeed, we would have that
\begin{align}
    \liminf_{t \to +\infty}\frac{2 \ma_H(\partial \Omega_t)}{( \frac{1}{16\pi}\int_{\partial \Omega_t} \H^2 \dif \sigma)^{\frac{1}{2}}+\frac{1}{16\pi}\int_{\partial \Omega_t} \H^2 \dif \sigma} = \liminf_{n \to +\infty} \ma_H(\partial \Omega_{t_n}) = \lim_{t \to +\infty} \ma_H(\partial \Omega_t),
\end{align}
where the last identity follows by \cref{thm:Geroch_monotonicity_formula}.

\medskip

To prove the claim, let $(t_n)_{n \in \N}$ be an increasing sequence divergent to $+\infty$ which realises the limit inferior in the right-hand side of \cref{eq:second_IMCF_de_lhopital}. By \cref{thm:Geroch_monotonicity_formula} we have two possible cases:
\begin{enumerate}
    \item there exists $T>0$ such that $\ma_{H}(\partial \Omega_{t_n})\geq 0$ for all $t_n \geq T$, or
    \item $\ma_{H}(\partial \Omega_{t_n})< 0$ for all $n\in \N$.
\end{enumerate}

\case Since $\ma_{H}(\partial \Omega_{t_n})\geq 0$, $\int_{\partial \Omega_{t_n}} \H^2 \dif \sigma\leq 16 \pi$ for every $t_n \geq T$.  By contradiction, suppose there exists $\varepsilon>0$ such that $\int_{\partial \Omega_{t_n}} \H^2 \dif \sigma \leq 16 \pi - \varepsilon$ for every $n$ sufficiently large. Then, by \cref{eq:second_IMCF_de_lhopital}, there exists $\kst(\varepsilon)>0$ such that
\begin{equation}
    +\infty>\liminf_{t\to +\infty} \ma_{\iso}(\Omega_{t}) \geq \lim_{n \to +\infty} \kst(\varepsilon)\sqrt{\abs{\partial \Omega_{t_n}}},
\end{equation}
which is clearly a contradiction. Hence, up to a not relabeled subsequence, $\int_{\partial \Omega_{t_n}} \H^2 \dif \sigma \to 16\pi$ as $n \to +\infty$. This proves the claim in this case.

\case Since $\ma_{H}(\partial \Omega_{t_n}) <0 $, $\int_{\partial \Omega_{t_n}} \H^2 \dif \sigma\geq 16 \pi$  for every $n\in \N$. Suppose by contradiction $\int_{\partial \Omega_{t_n}} \H^2 \dif \sigma\geq 16 \pi + \varepsilon$ for some $\varepsilon>0$. Then, by \cref{thm:Geroch_monotonicity_formula}, there exists $\kst(\varepsilon)>0$ such that
\begin{equation}
    \ma_{H}(\partial \Omega) \leq \lim_{t \to +\infty} \ma_{H}(\partial \Omega_t) \leq -\kst(\varepsilon)\lim_{n \to +\infty} \sqrt{\abs{ \partial \Omega_{t_n}}}=-\infty,
\end{equation}
which is a contradiction since $\h\in L^2(\partial \Omega)$, proving the claim also in this case.
\end{proof}

\begin{remark}
   In \cite[Appendix C]{chodosh_isoperimetryscalarcurvaturemass_2021}, the authors use an approach based on de l'H\^opital's rule to produce the following asymptotic comparison 
      \begin{equation}
      \label{eq:comparison-eichmair}
       \limsup_{V \to +\infty} \ma_{\iso}(E_V) \leq \limsup_{V\to +\infty} \ma_{H}(\partial E_V),
   \end{equation}
    where $E_V$ is the isoperimetric region of volume $V$. Their procedure uses the properties of the isoperimetric profile and sharp estimates on the Hawking mass, which hold at $\partial E_V$ provided it is connected.
    
  When isoperimetric sets have connected boundaries, one could combine \cref{eq:comparison-eichmair} with \cref{thm:mass_control_IMCF}  to deduce that both the Isoperimetric mass and the Hawking mass converge to  $\ma_{\iso}$ along the isoperimetric exhaustion as $V\to +\infty$. As a consequence, we would dispose of a different and simpler proof of \cref{thm:mass_equivalence} which does not use \cref{thm:bound_hawking_to_bound_iso} from \cite{jauregui_lowersemicontinuitymassconvergence_2019}. 
    
    To the authors' knowledge, this path is not feasible at the current state of the art. Indeed, it is unclear whether isoperimetric regions have connected boundaries under the sharp asymptotic assumptions we required on the metric $g$ in the statement of \cref{thm:mass_equivalence}. We refer the interested reader to \cite{benatti_isoperimetricsetsnonnegativescalar_2023} for more details.
   \end{remark}
   
    \begin{remark} When $\ma_H (\partial \Omega) \geq 0$, \cref{eq:main_inequality_IMCF} implies that the sublevel sets $\Omega_t$ asymptotically tend to satisfy a reverse sharp isoperimetric inequality. In fact, when the initial $\Omega$ is collapsed to a point $o \in M$, the techniques of Shi \cite{shi_isoperimetricinequalityasymptoticallyflat_2016} (after important insights by Brendle-Chodosh \cite{brendle_volumecomparisontheoremasymptotically_2014}) allow to show that the isoperimetric difference satisfies \begin{equation}
    \label{eq:shi}
    \abs{\Omega_t} -\frac{\abs{\partial \Omega_t}^{\frac{3}{2}}}{6 \sqrt{\pi}} \geq 0
    \end{equation}
    for the sublevel sets $\Omega_t$ of the weak IMCF issuing from $o\in M$ (we refer the reader to \cite{benatti_isoperimetricsetsnonnegativescalar_2023} for an explicit derivation). In this sense, \cref{eq:main_inequality_IMCF} can also be viewed as an asymptotic version of Shi's reverse isoperimetric inequality for weak IMCF starting at suitable initial surfaces. 
    \end{remark}

\subsection{Proof of the Isoperimetric Penrose Inequality and of the Isoperimetric PMT}
\label{sec:isoperimetricandtopological}
To finally prove \cref{thm:isoperimetric_penrose_intro}, we provide the following topological description  $3$-manifolds with nonnegative scalar curvature. It slightly improves on the analysis carried out in \cite[Section 4]{huisken_inversemeancurvatureflow_2001}, since we just assume the existence of an exhaustion of bounded sets with nonminimal mean-convex smooth boundaries, rather than the Asymptotically Flat condition. This assumption naturally fits the case of strongly $1$-nonparabolic manifolds, see \cref{rem:topology1nonpar}. The argument exploits the Mean Curvature Flow with surgery in Riemannian manifolds described by Brendle-Huisken \cite{brendle_meancurvatureflowsurgery_2018}, after \cite{brendle_meancurvatureflowsurgery_2016} and \cite{huisken_meancurvatureflowsingularities_1999}, and it was suggested in \cite[Remark, p. 394]{huisken_inversemeancurvatureflow_2001}. 
\begin{lemma}
\label{thm:topology}
Let $(M, g)$ be a Riemannian $3$-manifold with nonnegative scalar curvature and possibly with closed, smooth and minimal boundary. Assume that there exists an exhaustion $\{F_j\}_{j \in \N}$ of $M$ such that $\partial F_j$ are nonminimal, mean-convex closed surfaces.
\begin{enumerate}
    \item\label{item:outermost} If $\partial M$ is outermost, that is, there exists no other closed minimal surface homologous to $\partial M$, then we have $H_2(M,\partial M; \Z)=\set{0}$.
    \item\label{item:emptyboundary}  Assume that no minimal surface exists outside some compact set $K \subset M$, and that the closed surfaces $\partial F_j$ are spheres. Then there exists a possibly empty bounded $\Omega \subset M$ with smooth boundary $\partial \Omega$ such that  $M' = M \smallsetminus\Omega$ is a manifold with outermost minimal boundary $\partial M'=\partial \Omega$. As a consequence,  $H_2(M',\partial M' ;  \Z)=\set{0}$.
\end{enumerate}
\end{lemma}
\begin{proof}
First, we show that if $\partial M$ is outermost, then there are no closed surfaces with nonpositive mean curvature in $M$ other than $\partial M$. By contradiction, let $N$ be any such surface in $M$. Let $S$ be a nonminimal mean-convex smooth surface enclosing $N$. 
Then, we run the mean curvature flow starting at $S$ with surgery, as described by Brendle-Huisken \cite{brendle_meancurvatureflowsurgery_2018}. This procedure consists of carefully discarding necks of high mean curvature, which is the only type of singularity that arises in dimension $3$ other than spherical ones. We get an evolution of strictly mean-convex surfaces that are smooth except at finitely many surgery times. By \cite[Theorem 1.2]{brendle_meancurvatureflowsurgery_2018}, the flow will shrink to a closed minimal surface homologous to $\partial M$. Since $\partial M$ is outermost, the flow converges to $\partial M$. In particular, possibly delaying a surgery time, the mean-convex evolving sets would eventually touch $N$, giving a contradiction. 

In the assumption \cref{item:emptyboundary}, let $S$ be a nonminimal mean-convex sphere enclosing $K$. The mean curvature flow starting at $S$ ends in a closed minimal surface $\Sigma$ homologous to $\partial M$ or possibly becomes extinct at a finite collection of points. With the same argument as above, the region swept out by the mean curvature flow of $S$ is free of closed surfaces of nonpositive mean curvature. In particular, the manifold $M'$ has nonempty outermost minimal boundary $\Sigma$ or is a complete manifold without boundary, free of closed minimal surfaces.

For the rest of the proof, we focus on $M' = M \smallsetminus \Omega$. We agree that $M' =M$ if $\partial M$ is outermost minimal. We also assume that $\Sigma = \partial \Omega$ is nonempty since $\Sigma = \varnothing$ is a straightforward simplification. The Scalar Curvature Splitting Theorem \cite[Theorem 1.5]{galloway_scalarcurvaturewarpedproduct_2020} implies that $\Sigma$ must admit a metric with positive scalar curvature. Indeed, if such metric did not exist, $(M\smallsetminus {\Omega}, g)$ would isometrically split as a cylinder foliated by isometric copies of $\Sigma$. In that case, the cross-sections would be minimal, which contradicts the just proved absence of closed minimal surfaces. 

Since $\Sigma$ admits a metric with positive scalar curvature, each connected component of $\Sigma$ is a topological sphere, by the Gauss-Bonnet Theorem. Thus, we can run a slight variation of \cite[p.140]{lee_geometricrelativity_2019} to infer that $F_j$ satisfies $H_2(F_j \smallsetminus \Omega,\Sigma; \Z) = \set{0}$ for any $j$, which proves the lemma. Namely, consider the exact sequence
\begin{equation}
\label{sequence-homology}
    H_2(F_j \smallsetminus \Omega; \Z) \to H_2(F_j \smallsetminus \Omega, \Sigma; \Z) \to H_1(\Sigma ; \Z).
\end{equation}
By the mean-convexity of $\partial F_j$, we have that each nonzero class of $H_2(F_j \smallsetminus \Omega; \Z)$ is represented by a closed minimal surface (see e.g. \cite[Corollary 4.3]{lee_geometricrelativity_2019}). However, we just proved that $\Sigma$ is the only closed minimal surface, so the image of any class of $H_2(F_j \smallsetminus \Omega; \Z)$ vanishes in $H_2(F_j \smallsetminus\Omega, \Sigma; \Z)$. Moreover, since $\Sigma$ is a union of spheres, we have that $H_1(\Sigma ; \Z) = \{0\}$. The exact sequence \cref{sequence-homology} thus implies that $H_2(F_j \smallsetminus\Omega, \Sigma; \Z) =\{0\}$, completing the proof.
\end{proof}

\begin{remark}[Application to strongly $1$-nonparabolic $3$-manifolds]
\label{rem:topology1nonpar}
We observe that the assumptions in \cref{thm:topology}\cref{item:outermost} are satisfied in a strongly $1$-nonparabolic $3$-manifolds with nonnegative scalar curvature and closed outermost minimal boundary.  Indeed, by \cite[Minimizing Hull Property 1.4]{huisken_inversemeancurvatureflow_2001}, the sublevel sets $\{w \leq t\}$ of the weak IMCF starting at $\partial M$ are strictly outward minimising with $\CS^{1, \alpha}$ boundary, and they exhaust the manifold as $t \to + \infty$. We can approximate any of such boundaries with smooth strictly mean-convex surfaces by \cite[Lemma 2.6]{huisken_higherregularityinversemean_2008}, obtaining an exhaustion as in the statement.  
\end{remark}

\begin{remark} 
In \cref{thm:topology}, assume in addition that $M \smallsetminus K$ is diffeomorphic to $\R^3 \smallsetminus B$, for a suitable compact set $K$ and a ball $B$. Then, one can show that $M$ is diffeomorphic to $\R^3$ with a finite number of balls removed. This result can be deduced from the proof of \cite[Theorem 5.1]{eichmair_topologicalcensorship_2013} (cf. \cite[Theorem 4.11]{lee_geometricrelativity_2019}), which employs the Thurston's Geometrization of $3$-manifolds proved by Perelman \cite{perelman_entropyformularicciflow_2002, perelman_ricciflowsurgerythreemanifolds_2003}.
\end{remark}

We state an Isoperimetric version of the Positive Mass Theorem, which in particular does not require the boundary of the manifold to be outermost. 

\begin{theorem}[Isoperimetric Positive Mass Theorem]
\label{thm:isoperimetric_pmt}
Let $(M, g)$ be a complete strongly $1$-nonparabolic Riemannian $3$-manifold with nonnegative scalar curvature and possibly with closed, minimal boundary, such that the assumptions of \cref{thm:topology}\cref{item:emptyboundary} are met. Then,
\begin{equation}
\label{eq:isoperimetric_pmt}
   0 \leq \ma_{\iso},
\end{equation}
with equality achieved if and only if $(M, g)$ is isometric to the flat $\R^3$.
\end{theorem}
\begin{proof}[Proof of \cref{thm:isoperimetric_penrose_intro} and of \cref{thm:isoperimetric_pmt}]
 We first prove \cref{thm:isoperimetric_penrose_intro}. Let $N$ be one of the connected components of $\partial M$. By \cref{prop:jumping_to_homologous} there exists $\Omega$ with connected boundary homologous to $\partial M$ such that $\ma_{H}(\partial \Omega) \geq \ma_{H}(N)=\sqrt{\abs{ N}/16 \pi} > 0$ and $w_1\in \Lip_{\loc}(M \smallsetminus \Int \Omega )$ solution to \cref{eq:IMCF} starting at $\Omega$. \cref{thm:mass_control_IMCF,thm:Geroch_monotonicity_formula}, whose assumptions are fulfilled by \cref{rem:topology1nonpar} and \cref{thm:topology}\cref{item:outermost}, yield
\begin{equation}
    \sqrt{ \frac{\abs{N}}{16 \pi}} \leq \ma_H(\partial \Omega) \leq \lim_{t \to +\infty} \ma_H(\partial \Omega_t) \leq \liminf_{t \to +\infty} \ma_{\iso}(\Omega_{t}) \leq \ma_{\iso},
\end{equation}
where $\Omega_t = \set{ w_1 \leq t}$. This concludes the proof of \cref{eq:isoperimetric_penrose}.

We now prove \cref{eq:isoperimetric_pmt}, which, although weaker, does not require $\partial M$ to be outermost. 
By \cref{thm:topology}, \cref{thm:isoperimetric_penrose_intro} applies in the manifold $(M', g)$ with boundary $\partial M'$. Since the Isoperimetric mass of $(M',g)$ is the same as that of $(M, g)$, we get
\begin{equation}
\label{eq:penroseM'}
 \sqrt{\frac{ \abs{N'}}{16 \pi}} \leq \ma_{\iso},
\end{equation}
where $N'$ is a connected component of $\partial M'$. If $\partial M$ and $\partial M'$ are empty, take $o\in M$ and $w_1: M \smallsetminus \set{o} \to \R$ emanating from $o$. By \cref{thm:topology},  \cref{thm:Geroch_monotonicity_formula} applies. Hence, the Hawking mass is monotone along the weak IMCF starting from any initial set with connected boundary. Together with the construction considered in \cite[Lemma 8.1]{huisken_inversemeancurvatureflow_2001}, it yields $\ma_H(\partial \Omega_t)\geq 0$, where $\Omega_t = \set{w_1 \leq t}$. We conclude by \cref{thm:mass_control_IMCF} that $\ma_{\iso} \geq 0$ also in this case.

\smallskip

The rigidity statements in both \cref{thm:isoperimetric_penrose_intro} and \cref{thm:isoperimetric_pmt} stem from the constancy of the Hawking mass with the very same proof as in \cite[Proof of Main Theorem, step 2]{huisken_inversemeancurvatureflow_2001}. 
\end{proof}

\section{Application to Asymptotically Schwarzschildian manifolds}
\label{sec:asysvar}
As mentioned in the Introduction, the $\ADM$ concept of mass is well-defined only provided stronger asymptotic assumptions on the metric are fulfilled. If $(M,g)$ is the Schwarzschild of mass $\ma\geq 0$, the $\ADM$ mass coincide with $\ma$. As a consequence of the characterisation of isoperimetric sets obtained by Bray \cite{bray_penroseinequalitygeneralrelativity_1997}, it is not difficult to prove that also the Isoperimetric mass coincides with $\ma$. Since the volume and the perimeter of a set depend only on the $\CS^0$-character of the metric, it is natural to ask whether a Penrose inequality holds for $\ma$ when the space is merely $\CS^0$-asymptotic to the Schwarzschild of mass $\ma$. A further question is whether $\ma$ coincides with the Isoperimetric mass in this case.

\subsection{The weak IMCF on Asymptotically Flat Riemannian manifolds}
We introduce the concept of Asymptotically Flat Riemannian manifolds. This condition is of crucial importance in Mathematical General Relativity. Indeed, the asymptotically flat condition denotes the physical property that the gravitational system is isolated (see e.g. \cite[Section 3.1.2]{lee_geometricrelativity_2019}).

\begin{definition}[Asymptotically Flat Riemannian manifolds]
\label{def:asyflat}
A Riemannian $3$-manifold $(M,g)$ possibly with boundary is said to be \emph{$\CS^{k}_\tau$-Asymptotically Flat}, $k\in \N$ and $\tau >0$ ($\tau=0$ resp.) if the following conditions are satisfied.
\begin{enumerate}
    \item There exists a compact set $K \subseteq M$ such that $M \smallsetminus K$ is differmorphic to $\R^3\smallsetminus \set{\abs{x}\leq R}$, through a map $(x^1,x^2,x^3)$ whose component are called \emph{asymptotically flat coordinates}.
    \item In the chart $(M \smallsetminus K, (x^1,x^2,x^3))$ the metric tensor is expressed as
    \begin{equation}
        g= g_{ij} \dd x^i \otimes\dd x^j= (\eum_{ij}+\eta_{ij}) \dd x^i \otimes\dd x^j
    \end{equation}
    with 
    \begin{align}
      \qquad \sum_{i,j=1}^{3}\sum_{\abs{\beta}=0}^k \abs{x}^{\abs{\beta}+\tau} \abs{ \partial_\beta \eta_{ij}} =O(1) \text{ ($=o(1)$ resp.)} && \text{as } \abs{x} \to +\infty.
    \end{align}
\end{enumerate}
\end{definition}
Given an Asymptotically Flat Riemannian manifold $(M,g)$, $\D$ will denote the covariant derivative with respect to the metric $g$ and $\partial$ the one (which coincides with partial derivative) of $\eum$. We recall that $\abs{x}$ is asymptotically equivalent to the distance on $M$ as proved in \cite[Lemma 2.15]{benatti_asymptoticbehaviourcapacitarypotentials_2022}.

In a $\CS^0$-Asymptotically Flat Riemannian manifold, one can also compute the asymptotic expansion of the inverse of the metric in terms of the flat metric. Since $\eum^{ij}\eum_{jk}=g^{ij}g_{jk}=\delta^{i}_k$, we obtain that
\begin{equation}\label{eq:inverse_metrics}
    g^{il} - \eum^{il} = - \eum^{ij}\eta_{jk}g^{kl} = - \eum^{ij}\eta_{jk}\eum^{kl} +O(\abs{\eta}^2).
\end{equation}
In a $\CS^{1}$-Asymptotically Flat Riemannian manifold, observe that the Christoffel symbols of $g$ vanish at infinity as
\begin{equation}\label{eq:Christoffi}
    \Gamma^{k}_{ij} = \frac{1}{2} g^{kl} (\partial_i g_{lj} + \partial_j g_{li}- \partial_l g_{ij})= \frac{1}{2} g^{kl} (\partial_i \eta_{lj} + \partial_j \eta_{li}- \partial_l \eta_{ij}) = O(\abs{\partial \eta}).
\end{equation}

\medskip

$\CS^0$-Asymptotically Flat Riemannian manifolds are examples of strongly $1$-nonparabolic manifolds. Together with this assertion, the following proposition states an optimal control on the level sets. We will employ this tool in the proof of \cref{thm:penrose_adm_intro,thm:asysvarchild}.

\begin{proposition}
\label{prop:c01nonpar}
Let $(M, g)$ be a complete $\CS^0$-Asymptotically Flat Riemannian $3$-manifold. Then, $(M,g)$ is strongly $1$-nonparabolic. Moreover, consider any bounded $\Omega \subseteq  M$ with $\CS^1$-boundary homologous to $\partial M$ satisfying $\h \in L^2(\partial \Omega)$ and a solution $w_1$ to \cref{eq:IMCF}. Then, there exist $ \kst_1, \kst_2>0$ such that 
\begin{equation}
\label{eq:boundimcf}
2\log \abs{x} - \kst_1 \leq w_1(x) \leq 2 \log \abs{x} + \kst_2   
\end{equation}
for any $x \in M \smallsetminus K$, where $K$ is the compact set of \cref{def:asyflat}.
\end{proposition}
\begin{proof}
    We deferred the discussion of the strongly $1$-nonparabolic nature of $(M,g)$ to \cref{rem:C01pnonpar}. Given $\Omega$ with smooth boundary, \cref{thm:stronglyp-li-yau} shows that solutions $(w_p)_{p>1}$ of \cref{eq:pIMCF} satisfy an uniform lower bound that diverges as $\abs{x} \to +\infty$. Hence, the locally uniform limit $w_1$ obtained 
    as in \cite{moser_inversemeancurvatureflow_2007, moser_inversemeancurvatureflow_2008, kotschwar_localgradientestimatesharmonic_2009, mari_flowlaplaceapproximationnew_2022} sending $p\to1^+$ is the unique proper IMCF starting at $\Omega$. The double bound \cref{eq:boundimcf} then follows from \cref{eq:li-yau-p}. For $\Omega$ with $\CS^1$-boundary, one can approximate $\partial \Omega$ in $\CS^1$ with smooth surfaces in $M \smallsetminus \partial \Omega$ and use the solutions starting at those surfaces as barriers.
\end{proof}

In the following remark, we specialize \cref{thm:topology} to describe the topology of asymptotically flat $3$-manifolds.

\begin{remark}[Topology of asymptotically flat manifolds]
\label{rem:C1-foliation}
Firstly, \cref{prop:c01nonpar} grants that $\CS^0$-Asymptotically Flat $3$-manifolds are strongly $1$-nonparabolic. Hence, any $\CS^0$-Asymptotically Flat $3$-manifold with nonnegative scalar curvature and outermost minimal boundary satisfies the assumptions of \cref{thm:topology}\cref{item:outermost} (see \cref{rem:topology1nonpar}). This fact will be employed in the proof of \cref{thm:asysvarchild,thm:penrose_adm_intro}. 

\smallskip

Secondly, the assumptions in \cref{thm:topology}\cref{item:emptyboundary} are fulfilled by any $\CS^1$-Asymptotically Flat manifold with nonnegative scalar curvature and outermost minimal boundary. In this case, the level sets of $\abs{x}$ constitute a strictly mean-convex foliation at infinity, whose presence prevents the existence of closed minimal surfaces in the foliated region. Indeed, the divergence operator only involves the metric $g$ up to the first derivatives. Hence, \cref{eq:Christoffi} ensures
    \begin{equation}
    \label{eq:strictmeanconv}
        \div \left(\frac{\D \abs{x}}{\abs{\D \abs{x}}}\right) > 0
    \end{equation}
    for $\abs{x}$ sufficiently large. 
    We will use this property to prove the Positive Mass Theorem for the $\ADM$ mass \cref{thm:pmt_adm}. 
\end{remark}

\subsection{Riemannian Penrose Inequality on Asymptotically Schwarzschildian manifolds}

We now provide the precise definition of Asymptotically Schwarzschildian manifolds, which is the setting of \cref{thm:asysvarchild}. As for the Asymptotically Flat condition, one could give the definition of $\CS^k_\tau$-Asymptotically Schwarzschildian manifolds for every $k\in \N$, but since we will never use this general notion, we prefer to focus on the case $k=0$.

\begin{definition}[Asymptotically Schwarzschildian manifolds]
\label{def:asysvar}
A Riemannian $3$-manifold $(M,g)$ possibly with boundary is \emph{$\CS^{0}_\tau$-Asymptotically Schwarzchildian of mass $\ma\geq 0$}, $\tau >0$ ($\tau=0$ resp.), if the following conditions are satisfied.
\begin{enumerate}
    \item There exists a compact set $K \subseteq M$ such that $M \smallsetminus K$ is diffeomorphic to $\R^3\smallsetminus \set{\abs{x}\leq R}$, through a map $(x^1,x^2,x^3)$ whose component are called \emph{asymptotically flat coordinates}.
    \item In the chart $(M \smallsetminus K, (x^1,x^2,x^3))$ the metric tensor is expressed as
    \begin{equation}    \label{eq:asy-svar}
        g= g_{ij} \dd x^i \otimes\dd x^j= \left(\left(1+ \frac{2\ma}{\abs{x}}\right)\eum_{ij}+\eta_{ij}\right) \dd x^i \otimes\dd x^j
    \end{equation}
    with
    \begin{align}
       \sum_{i,j=1}^{3}\abs{x}^{1+\tau} \abs{ \eta_{ij}} =O(1) \text{ ($=o(1)$ resp.)} && \text{as $\abs{x} \to +\infty$}.
    \end{align}
\end{enumerate}
\end{definition}

Observe that the metric $g$ can be expressed as
\begin{align}
    g_{ij}=\left(1+ \frac{\ma}{2\abs{x}}\right)^4\eum_{ij}+O(\abs{x}^{-\min\set{\tau+1,2}}).
\end{align}
Then, the metric $g$ is Asymptotically equivalent to the Schwarzschild metric at infinity, at least to the first order. Moreover, observe that a $\CS^0_\tau$-Asymptotically Schwarzschildian manifold is always $\CS^0_1$-Asymptotically Flat.

\medskip

We are finally ready to prove the Penrose inequality for $\CS^0_\tau$-Asymptotically Schwarzschildian manifolds. 

\begin{proof}[Proof of \cref{thm:asysvarchild}]
We assume without loss of generality that $\tau \leq 1$.
It is easy to show that 
\begin{equation}
\label{eq:limpalle}
\lim_{R \to + \infty}\ma_{\iso}(\set{\abs{x} \leq R}) = \ma,
\end{equation}
which directly implies
\begin{equation}
    \label{eq:mminoremiso}
    \ma \leq \ma_{\iso}.
\end{equation}
We briefly provide the computations involved. Given $\Sigma = \partial \Omega$ a closed surface in $M$, the measure $\dif \sigma$ induced by $g$ on $\Sigma$ is related to the one induced by $\eum$ by
\begin{equation}
    \label{eq:sviluppoaree}
    \dif\sigma= \left(1 + \frac{2\ma}{\abs{x}}\right)(1 + \eum_{\top}^{ij}\eta_{ij} + O(\abs{\eta}^2))\dif\sigma_\eum, 
\end{equation}
where $\eum_{\top}^{ij} = \eum^{ij} - \nu^i_{\eum}\nu^j_{\eum}$ with $\nu_{\eum}$ outer $\eum$-unit normal to $\Sigma$.
Applying this formula, we immediately get 
\begin{equation}
    \label{eq:asyarea}
    \abs{\set{\abs{x} = R}} =  4\pi R^2 \left(1 + \frac{2\ma}{R} + O(\abs{R}^{-1-\tau}) \right).
\end{equation}
Moreover, we have
\begin{equation}
    \label{eq_asygrad}
    \abs{\D \abs{x}}^2= 1 - \frac{2\ma}{\abs{x}} + O(\abs{x}^{-1-\tau}),
\end{equation}
from which
\begin{equation}
    \label{eq:asyinversegrad}
    \abs{\D \abs{x}}^{-1} = 1 + \frac{\ma}{\abs{x}} + O(\abs{x}^{-1-\tau})
\end{equation}
follows. Then, coupling with \cref{eq:asyarea}, it holds
\begin{equation}
\label{eq:derivatavolume}
    \frac{\dif}{\dif R} \abs{\set{\abs{x} \leq R}} = \int\limits_{\set{\abs{x} = R}} \frac{1}{\abs{\D \abs{x}}} \dif \sigma= 4\pi R^2 + 12 \pi \ma R + O(R^{1-\tau}).
\end{equation}
Hence, by integrating, we obtain
\begin{equation}
    \label{eq:asyvolume}
    \abs{\set{\abs{x} \leq R}} = \frac{4}{3} \pi R^3 + 6\pi \ma R^2 + O(R^{2-\tau}).
\end{equation}
Finally, by \cref{eq:asyvolume} and \cref{eq:asyarea} we compute 
\begin{equation}
\label{eq:finalasymballs}
\begin{split}
\ma_{\iso}(\set{\abs{x} \leq R}) =  \frac{2}{\abs{\set{\abs{x} = R}}} \left(\abs{\set{{\abs{x} \leq R}}} - \frac{\abs{\set{\abs{x} = R}}^\frac{3}{2}}{6\sqrt{\pi}}  \right) 
= \frac{\ma (1 + O(R^{-\tau}))}{1 + O(R^{-1})}.
\end{split}
\end{equation}
Letting $R \to +\infty$ in \cref{eq:finalasymballs} yields \cref{eq:limpalle}.

We are now left to prove that $\ma_{H}(\partial \Omega) \leq \ma$ for every $\Omega\subseteq M$ bounded with connected $\CS^1$-boundary homologous to $\partial M$ satisfying $\h\in L^2(\partial \Omega)$. Indeed, by \cref{prop:jumping_to_homologous,thm:bound_hawking_to_bound_iso} this implies $\ma_{\iso}\leq \ma$ that coupled with \cref{eq:mminoremiso} gives $\ma=\ma_{\iso}$. The Riemannian Penrose Inequality now follows by \cref{thm:isoperimetric_penrose_intro}.

We consider the weak IMCF $w_1$ emanating from some bounded $\Omega$ with smooth connected $\CS^1$-boundary homologous to $M$. We let $\Omega_t =\set{w_1\leq t}$. 
We denote with  $\svm$ the Schwarzschild metric, that from its very definition satisfies
\begin{equation}
    \label{eq:asysvar}
    \svm_{ij} = \left(1 + \frac{2\ma }{\abs{x}}\right)\eum_{ij}+ O(\abs{x}^{-2})
\end{equation}
for every $i,j=1,2,3$ as $\abs{x} \to +\infty$. Then, applying again \cref{eq:sviluppoaree}, we immediately get
\begin{equation}
    \label{eq:comparisonaree-svar}
    \abs{\abs{\partial \Omega_t} - \abs{\partial \Omega_t}_{\svm}} =\abs{ \int_{\partial \Omega_t} O(\abs{x}^{-1-\tau} )\dif\sigma_\eum}.
\end{equation}
By the second inequality in \cref{eq:boundimcf} we have $\abs{x}^{2} \geq \kst \ee^t$ for any $x \in \partial \Omega_t$, for $t$ large enough. Thus \cref{eq:comparisonaree-svar} implies
\begin{equation}
    \label{eq:asyareesvarpromoted}
     \abs{\abs{\partial \Omega_t} - \abs{\partial \Omega_t}_{\svm}} \leq \kst \abs{\partial \Omega_t}_{\eum} \ee^{-t \frac{(1+ \tau)}{2}}. 
     \end{equation}
     The very same argument also shows that
     \begin{equation}
         \label{eq:asyR3svar}
         \abs{\abs{\partial \Omega_t}_{\svm} - \abs{\partial \Omega_t}_{\eum}} \leq \kst\abs{\partial \Omega_t}_{\eum} \ee^{-t}.
        \end{equation}
As a consequence, coupling \cref{eq:asyareesvarpromoted} and \cref{eq:asyR3svar}, we can write
\begin{equation}
    \label{eq:misomodified}
    \ma_{\iso}(\Omega_t) = \frac{2}{\abs{\partial \Omega_t}}\left(\abs{\Omega_t}- \frac{\abs{\partial \Omega_t}^{\frac{3}{2}}_{\svm}}{6\sqrt{\pi}} + O\left(\abs{\partial \Omega_t}_{\svm}^{\frac{3}{2}} \ee^{-t \frac{(1+\tau)}{2}}\right) \right).
\end{equation}
We are going to pass to the superior limit as $t \to +\infty$ in \cref{eq:misomodified}. Let $B_t$ be the coordinate ball in the Schwarzschild space of volume $\abs{B_t}_{\svm} = \abs{\Omega_t}$. Since $B_t$ is isoperimetric in the Schwarzschild metric for its own volume \cite[Theorem 8]{bray_penroseinequalitygeneralrelativity_1997} \cite[Corollary 2.6]{bray_isoperimetriccomparisontheoremschwarzschild_2002}, we have
\begin{equation}
    \label{eq_limpart1}
    \begin{split}
    \frac{2}{\abs{\partial \Omega_t}}\left(\abs{\Omega_t} - \frac{\abs{\partial \Omega_t}^{\frac{3}{2}}_{\svm}}{6\sqrt{\pi}}\right) &= \frac{2\abs{\partial \Omega_t}_{\svm}}{\abs{\partial \Omega_t}\abs{\partial \Omega_t}_{\svm}}\left(\abs{B_t}_{\svm} - \frac{\abs{\partial \Omega_t}^{\frac{3}{2}}_{\svm}}{6\sqrt{\pi}}\right) \\ 
    &\leq \frac{2\abs{\partial \Omega_t}_{\svm}}{\abs{\partial \Omega_t}\abs{\partial B_t}_{\svm}}\left(\abs{B_t}_{\svm} - \frac{\abs{\partial B_t}^{\frac{3}{2}}_{\svm}}{6\sqrt{\pi}}\right).
\end{split}
\end{equation}
Since $\abs{\partial \Omega_t}$ is asymptotic to $\abs{\partial \Omega_t}_{\svm}$ by \cref{eq:asyareesvarpromoted} and \cref{eq:asyR3svar}, \cref{eq:limpalle} yields  
\begin{equation}
\label{eq:limpart1ok}
    \limsup_{t \to + \infty} \frac{2}{\abs{\partial \Omega_t}}\left(\abs{\Omega_t} - \frac{{\abs{\partial \Omega_t}}_{\svm}^{\frac{3}{2}}}{6\sqrt{\pi}}\right) \leq \ma. 
\end{equation}
We are left to estimate the remaining part of \cref{eq:misomodified}. First, observe that 
\begin{equation}
    \label{limpart2}
    \ee^{-t\frac{(1 + \tau)}{2}}\frac{\abs{\partial \Omega_t}^{\frac{3}{2}} }{\abs{\partial \Omega_t}} = \ee^{-t\tau}\abs{\partial \Omega^*}^{\frac{1}{2}} 
 \end{equation}
 since $\abs{\partial \Omega_t}= \abs{\partial \Omega^*} \ee^t$. Hence, 
    \begin{equation}\label{limpart2fine}
      \lim_{t \to +\infty}  
       \ee^{-t\frac{(1 + \tau)}{2}} \frac{\abs{\partial \Omega_t}^{\frac{3}{2}} }{\abs{\partial \Omega_t}}= 0. 
    \end{equation}
Since $\abs{\partial \Omega_t}$ and $\abs{\partial \Omega_t}_{\svm}$ are asymptotically equivalent by \cref{eq:asyareesvarpromoted} and \cref{eq:asyR3svar}, plugging \cref{eq:limpart1ok} and \cref{limpart2fine} into \cref{eq:misomodified} infers
\begin{equation}
    \label{eq:misoninfinity}
    \limsup_{t \to + \infty} \ma_{\iso} (\Omega_t) \leq \ma.
\end{equation}
Thanks to  \cref{thm:topology} and \cref{thm:mass_control_IMCF}, \cref{thm:Geroch_monotonicity_formula} applies and concludes the proof. \end{proof}
It would be very nice to remove the request of nonnegative mass in the statement of \cref{thm:asysvarchild}. Such an assumption has been crucially exploited when estimating from above the quasi-local Isoperimetric mass of an evolving set $\Omega_t$ with that of a geodesic ball in the space-like Schwarzschild metric $s$ with the same volume. It was possible because such balls are isoperimetric in the Schwarzschild metric with nonnegative mass, while they are not if $\ma$ is negative, as shown in \cite[Section 5]{corvino_isoperimetricsurfacesgeneralrelativity_2007}. Understanding the shape of isoperimetric sets in the Schwarzschild metric with negative mass could be crucial in showing that, $\ma \geq 0$ is a consequence in \cref{thm:asysvarchild} rather than an assumption.

\section{Riemannian Penrose Inequality for the \texorpdfstring{$\ADM$}{ADM} mass}
\label{sec:admpenrose}
In this section, we prove the Riemannian Penrose Inequality \cref{thm:penrose_adm_intro} for the $\ADM$ mass, in the natural $\CS^{1}_\tau$-Asymptotically Flat setting, $\tau > 1 /2$. We take advantage of the fine analysis at infinity for harmonic functions to infer a finite bound for the Hawking masses of surfaces $\partial\Omega$ homologous to $\partial M$. Indeed, our proof of \cref{thm:penrose_adm_intro} combines the Geroch monotonicity \cref{thm:Geroch_monotonicity_formula} with its potential-theoretic counterpart. Such monotonicity was discovered in \cite{agostiniani_greenfunctionproofpositive_2024}. Here, the authors show that the \emph{$2$-Hawking mass}
\begin{equation}
\label{eq:2hawking}
 \ma^{\tph}_{H}(\partial \Omega) = \frac{\ncapa_2(\partial \Omega)}{8 \pi}\left[4\pi + \int_{\partial \Omega}{\abs{\D w_2}^2}\dif \sigma - \int_{\partial \Omega}{\abs{ \D w_2} }\H \dif \sigma \right]
\end{equation}
is monotone along the level sets of $w_2$, which solves
\begin{equation}\label{eq:m-potential}
    \begin{cases}
     \Delta w_2& =& \abs{\D w_2}^2 &\text{on $M \smallsetminus\overline{\Omega}$,}\\
     w_2&=&0 & \text{on $\overline{\Omega}$,}\\
     w_2&\to& +\infty & \text{as $\abs{x}\to +\infty$.}
    \end{cases}
\end{equation}

The existence of $w_2$ is granted by the Asymptotically Flat condition (see \cref{thm:stronglyp-li-yau}).
In what follows, we are exploiting the linear nature of the PDE solved by $u = - \log w_2$ (see \cref{eq:p-cap-potential}) to obtain a deeper asymptotic understanding of the quantities involved in \cref{eq:2hawking}. This approach circumvents the additional asymptotic assumption on the Ricci curvature tensor, as mandated by Huisken-Ilmanen \cite{huisken_inversemeancurvatureflow_2001} and Agostiniani-Mantegazza-Mazzieri-Oronzio \cite{agostiniani_riemannianpenroseinequalitynonlinear_2022}. In detail, we can control the $L^2$-norm of the Hessian of any function with the $L^2$-norms of its gradient and Laplacian by exploiting the $\CS^1$-asymptotic behaviour of the metric. Utilising this estimate alongside the representation formula, we obtain the expected decay of the error between the function $u=-\log w_2$ and its rotationally symmetric model on flat $\R^3$ up to its first derivatives. Then, we improve it to its second-order derivatives, hence to the mean curvature of its level sets, in an $L^2$-fashion.

Such knowledge is sufficient to bound $\ma^{\tph}_{H}(\partial \Omega_t)$ from above in terms of $\ma_{\ADM}$ (see \cref{thm:bound_2_hawking}) and, in turn, it yields a finite, uniform (nonsharp) bound for the Hawking mass of any initial set (see \cref{thm:hawking_p-hawking}). We refine such bound through an improved asymptotic behaviour of the curvatures of the sets evolving through weak IMCF. Here, it is where we make up for the milder asymptotic assumptions on the metric. Indeed, we can conclude 
\begin{equation}\label{eq:hawking_madm_intro}
 \ma_{H}(\partial \Omega) \leq \ma_{\ADM}
\end{equation}
for a relevant class of subsets of $M$ as in the final analysis of \cite[Asymptotic Comparison Lemma 7.4]{huisken_inversemeancurvatureflow_2001}, but 
without requiring that $(M,g)$ is $\CS^{1}_1$-Asymptotically Flat.

The inequality
\cref{eq:hawking_madm_intro} implies at once the Penrose Inequality \cref{eq:adm_penrose} when applied to $\partial \Omega =N $. Combining \cref{eq:hawking_madm_intro} with Jauregui-Lee's \cref{thm:bound_hawking_to_bound_iso} and with the limit behaviour of the quasi-local Isoperimetric masses of geodesic balls computed in \cite{fan_largespheresmallspherelimitsbrownyork_2009}, we also obtain \cref{thm:mass_equivalence}.

\subsection{Preliminaries on Linear Potential Theory} 
Let $(M,g)$ be a $\CS^0$-Asymptotically Flat Riemannian $3$-manifold. Consider a bounded subset $\Omega\subseteq M$ with $\CS^{1,\alpha}$-boundary homologous to $\partial M$ and with $\h \in L^2(\partial \Omega)$. Then, a solution $w_2\in \CS^\infty(M\smallsetminus \Omega)\cap\CS^{1,\alpha}_{\loc}(M\smallsetminus \Int \Omega)$ to \cref{eq:m-potential} exists. Indeed, $w_2= -\log u$ where $u$ is the capacitary potential of $\Omega$ and solves
\begin{equation}\label{eq:p-cap-potential}
    \begin{cases}
     \Delta u& =& 0 &\text{on $M \smallsetminus \overline{\Omega}$,}\\
     u&=&1 & \text{on $\overline{\Omega}$,}\\
     u&\to& 0 & \text{as $\abs{x}\to +\infty$.}
    \end{cases}
\end{equation}

The capacity is defined for any compact $K \subset M$ as the quantity
\begin{equation}\label{eq:p-capacity}
    \ncapa_2(K) = \inf \set{\frac{1}{4\pi}\int_M \abs{\D v}^2 \dif \mu\st v \in \CS_c^\infty(M),\, v \geq 1\text{ on } K}.
\end{equation}
It turns out that the capacity of the level sets of solutions to \cref{eq:m-potential} exponentially grows, exactly as the perimeter of level sets of the weak IMCF does. We recall this useful property in the following lemma (see \cite[Lemma 3.8]{holopainen_nonlinearpotentialtheoryquasiregular_1990}).

\begin{lemma}\label{lem:capacitygrowth}
Let $(M,g)$ be a Riemannian $3$-manifold possibly with boundary. Let $\Omega \subseteq M$ be bounded with boundary homologous to $\partial M$. Given a solution $w_2$ to \cref{eq:m-potential}, denote $\Omega_t = \set{w_2 \leq t}$. We have
\begin{equation}
    \ncapa_2(\partial \Omega_t) = \ee^{t} \ncapa_2(\partial \Omega)= \frac{1}{4\pi}\int_{\partial \Omega_t} { \abs{ \D w_2}}\dif \sigma.
\end{equation}
\end{lemma} 

The $2$-Hawking mass \cref{eq:2hawking} represents the potential-theoretic version of the Hawking mass, and it is monotone along the level set flow of the function $w_2$ as shown in \cite{agostiniani_greenfunctionproofpositive_2024}. More precisely, observe that the following relations are true
\begin{equation}\label{eq:literature_mass_relations}
    \frac{\ncapa_2(\partial \Omega)}{8\pi} U_2\left(\ee^{t}\right) = \ma_{H}^{\tph}(\partial \Omega_t)= \frac{1}{8\pi}F_2\left(\ncapa_2(\partial \Omega) \ee^{t}\right),
\end{equation}
where $F_2(t)$ is the function defined in \cite[(1.2)]{agostiniani_greenfunctionproofpositive_2024} and $U_2(t)$ is the function
\begin{equation}\label{eq:mass_potential}
    U_2(t)= 4 \pi t + t^{3} \smashoperator{\int\limits_{\set{u = 1/t}}} \abs{ \D u}^2 \dif \sigma -t^{2}  \smashoperator{\int\limits_{\set{u =1/t}}} \abs{ \D u}\H \dif \sigma.
\end{equation}
The following result is then a consequence of \cite[Theorem 1.1]{agostiniani_greenfunctionproofpositive_2024} (see also \cite{agostiniani_riemannianpenroseinequalitynonlinear_2022}).

\begin{theorem}[Monotonicity of the $2$-Hawking mass]
\label{thm:AMMO_monotonicity}
    Let $(M,g)$ be a complete $\CS^0$-Asymptoti\-cally Flat Riemannian $3$-manifold with nonnegative scalar curvature and possibly with smooth and closed boundary. Assume that $H_2(M, \partial M;\Z) = \set{0}$. Let $\Omega\subseteq M$ with connected $\CS^{1,\alpha}$-boundary homologous to $\partial M$ satisfying $\h \in L^2(\partial \Omega)$ and $w_2$ the solution to \cref{eq:m-potential} starting at $\Omega$. Then, denoting $\Omega_t = \set{w_2 \leq t}$, the function $t\mapsto \ma_H^{\tph}(\partial \Omega_t)$ defined in \cref{eq:2hawking} belongs to $\BV_{\loc}(0,+\infty)$ and admits a monotone nondecreasing representative. 
\end{theorem}

We do not claim that we can jump through horizons, since it is not clear to us whether a result analogous to \cite[Geroch Monotonicity (Multiple Boundary Components) 6.1]{huisken_inversemeancurvatureflow_2001} holds along the level set flow of $w_2$. On the other hand, the following elementary control with the classical Hawking mass will suffice for our aims.

\begin{lemma}\label{thm:hawking_p-hawking}
    Let $(M,g)$ be a complete $\CS^0$-Asymptotically Flat Riemannian $3$-manifold possibly with smooth and closed boundary. Then, for every outward minimising $\Omega\subseteq  M$ with $\CS^{1,\alpha}$-boundary homologous to $\partial M$ with $\h \in L^2(\partial \Omega)$ we have
    \begin{equation}
         \ma_{H}(\partial \Omega)\leq \kst\ma^{\tph}_{H}(\partial \Omega),
    \end{equation}
    where $\kst$ depends only on the isoperimetric constant of $(M,g)$.
\end{lemma}

\begin{proof} Observe that
\begin{align}
     \int_{\partial \Omega}\abs{\D w_2}\H  \dif \sigma - \int_{\partial \Omega}\abs{\D w_2}^2\dif \sigma =\int_{\partial \Omega} \frac{\H^2}{4} \dif \sigma - \int_{\partial \Omega} \left( \frac{\H}{2} - \abs{ \D w_2}\right)^2 \dif \sigma\leq  \int_{\partial \Omega} \frac{ \H^2}{4} \dif \sigma.
\end{align}
To conclude is then enough to proceed as in \cite[Theorem 1.3]{fogagnolo_minimisinghullscapacityisoperimetric_2022} to prove that
\begin{equation}
    \ncapa_2 (\partial \Omega) \geq \kst \abs{ \partial \Omega}^{\frac{1}{2}}.\qedhere
\end{equation}
\end{proof}

\subsection{Sharp integral asymptotic estimates for harmonic functions}
On $\CS^1_\tau$-Asymptotical\-ly Flat manifolds, it seems out of reach to obtain $\CS^2$-asymptotic information for the capacitary potential at infinity. However, we are actually deriving quantitative (in terms of $\tau$) $L^2$- asymptotic estimates for the Hessian at infinity. The strategy of the proof consists in a delicate refinement of the classical one, clearly exposed in \cite[Appendix A]{mantoulidis_capacityquasilocalmasssingular_2020} (see also \cite[Lemma 4.1]{hirsch_positivemasstheoremmanifolds_2020}).

The following direct consequence of Bochner's identity will be a basic tool. We exploit it using an approach partially similar to the one described in \cite{fogagnolo_notecriticallaplaceequation_2022}. We work for simplicity in $\R^3$ endowed with an Asymptotically Flat metric; this will simplify the notation and the presentation and, for our purposes, it will clearly result in no loss of generality.

\begin{lemma}
\label{lem:boch}
Let $k>1$, $R_0>0$ and $\varepsilon>0$ fixed. There exists a constant $\kst>0$, such that for any metric $g$ on $\R^3$ with $\set{\abs{x} \geq R_0}\subseteq \set{\abs{g-\delta} + \abs{x}\abs{\partial g}\leq \varepsilon } $ we have
\begin{equation}
    \label{eq:boch1}
    \int\limits_{A_{R, kR}} \abs{ \D\D f}^2 \dif \mu \leq \kst \left( \frac{1}{R^2}\int\limits_{A_{R/2, 2kR}} \abs{\D f}^2 \dif\mu+ \int\limits_{A_{R/2, 2kR}} (\Delta f)^2 \dif \mu \right),
    \end{equation}
for any $f \in \CS^2(\R^3)$ and all $R \geq 2 R_0$, where $A_{R,kR}=\set{R \leq \abs{x}\leq kR}$.
\end{lemma}

\begin{remark} 
Observe that if $g=\delta$, the above lemma holds for $R_0=0$. In the general case, up to choosing $R_0$ sufficiently large, \cref{eq:boch1} holds with the same constant regardless of whether the terms are expressed with respect to the metric $g$ or $\eum$. 
\end{remark}

\begin{proof}
Suppose that $g= \eum$. Let $\varphi$ be the classic cut-off function supported in $A = A_{R/2, 2k R}$ which is equal to $1$ on $A_{R, kR}$ with $\abs{\D \varphi} \leq \kst/R$ and $\abs{\Delta \varphi} \leq \kst/R^2$. Multiplying the Bochner identity
\begin{equation}
    \label{eq:boch-classical}
    \sum_{ij}(\partial_i \partial_j f)^2 = \sum_{ij}\frac{1}{2} \partial_i^2 (\partial_jf)^2 - \partial_j\partial_i^2f \partial_j f 
\end{equation}
by $\varphi$ and integrating on $A$, we get
\begin{equation}
    \label{eq:integralboch}
    \begin{split}
    \int_A \abs{\partial \partial f}_{\eum}^2 \varphi \dif x &\leq \int_A \abs{\partial f }_{\eum}^2 \abs{\Delta_\eum \psi} + (\Delta_\eum f)^2 \varphi + \abs{\partial f }_{\eum}\abs{\Delta_{\eum} f} \abs{\partial \varphi}_{\eum} \dif x\\
    &\leq\kst \left( \frac{1}{R^2}\int_{A} \abs{\partial f}_{\eum}^2 \dif x + \int_A (\Delta_{\eum} f)^2 \dif x \right) 
    \end{split}
\end{equation}
after performing integration by parts and Cauchy-Schwarz inequality. Hence, the theorem follows for $g = \eum$.
Let now $g$ be as in the statement. In particular, $g_{ij}=\eum_{ij} + \eta_{ij}$ with $\abs{\eta}_{\delta}+\abs{x}{\abs{\partial \eta}}_{\delta}\leq \varepsilon$ on $\set{\abs{x} \geq R_0}$. On $A$, we have 
\begin{equation}
\label{eq:hessian-close}
\begin{split}
 \abs{\D \D f}^2_g  &\leq \kst\abs{ \partial_i \partial_jf -\Gamma_{ij}^k \partial_kf}^2_\eum \leq \kst\left(\abs{ \partial \partial f}^2_{\eum} +\abs{\partial \eta}_{\eum}\abs{\partial \partial f}_{\eum} \abs{ \partial f}_\eum+ \abs{\partial \eta}^2_{\eum} \abs{\partial f}^2_{\eum}\right)\\
 &\leq \kst\left( \abs{\partial \partial f}_{\eum}^2 + \abs{\partial \eta}^2_{\eum} \abs{ \partial f}^2_{\eum}\right)\leq  \kst \abs{ \partial \partial f}^2_{\eum} +  \frac{\kst}{R^2}\abs{\partial f}_{\eum}^2,
 \end{split}
\end{equation}
where $\kst$ does not depends on $R$. Similarly,
\begin{align}
    \abs{\partial f}_{\eum}^2 \leq \kst \abs{ \D f }_g^2, && \abs{\Delta_{\eum} f} \leq \kst \abs{ \Delta_g f} + \frac{ \kst}{R} \abs{ \D f}_g
\end{align}
hold. Employing all the above estimates we get
\begin{align}
    \int_{A} \abs{ \D\D f}^2 \dif \mu &\leq \kst \left(\int_{A}  \abs{ \partial \partial f}^2_{\eum} + \frac{1}{R^2}\abs{\partial f}_{\eum}^2\dif \mu \right)\leq \kst \left(\int_{A}  \abs{ \partial \partial f}^2_{\eum} \dif x +\int_{A}  \frac{1}{R^2}\abs{\D f}_g^2\dif \mu \right)\\
    &\leq \kst \left(\int_{A}  (\Delta_{\eum} f)^2 + \frac{1}{R^2} \abs{ \partial f}^2 \dif x +\int_{A}  \frac{1}{R^2}\abs{\D f}_g^2\dif \mu \right)\\
    &\leq \kst \left(\int_{A}  (\Delta f)^2 \dif \mu +\int_{A}  \frac{1}{R^2}\abs{\D f}_g^2\dif \mu \right),
\end{align}
where we used the fact that $\kst^{-1}\leq \sqrt{ \det g} \leq \kst$ by a constant that does not depend on $R$. \end{proof}

\smallskip

In the sequel of this section, we fix $k > 1$ and we derive decay estimates for a solution $u$ to \cref{eq:p-cap-potential}. One might rephrase all the computations below directly in terms of solutions $w_2$ to \cref{eq:m-potential}. However, we think that the proposed presentation is simpler and the reader would more readily compare our results with the ones in the literature.

\begin{lemma}
\label{lem:decayL2}
Let $0 < \tau  \leq 1$ and let $(\R^3,g)$ be a $\CS^1_\tau$-Asymptotically Flat Riemannian manifold. Then, the capacitary potential $u$ \cref{eq:p-cap-potential} of a bounded set $\Omega$ with smooth boundary satisfies the  $\CS^1$-decay 
\begin{align}
    \label{eq:decay-pointwise}
        u(x)=O(\abs{x}^{-1})&&
        {\abs{\D u}(x)} = O(\abs{x}^{-2})
\end{align}
as $\abs{x} \to +\infty$. Consequently, the $L^2$-asymptotic estimate 
\begin{equation}
\label{eq:decayL2}
\int\limits_{A_{R, kR}} \abs{\D \D u}^2 \dif\mu =O(R^{-3})
\end{equation}
holds as $R\to +\infty$, where $A_{R,kR}=\set{R \leq \abs{x}\leq kR}$.
\end{lemma}
\begin{proof}
The $\CS^0$ decay of $u$ follows exactly as in the first part of the proof \cite[Lemma A.2]{mantoulidis_capacityquasilocalmasssingular_2020}. The $g$-Laplacian of $v = a \abs{x}^{-1} - \abs{x}^{-1- \varepsilon}$ (that depends only on the $\CS^0$-character of $g$) is nonpositive for $0 <\varepsilon < \tau$, and, as such, choosing $a$ big enough we get $u \leq v \leq \kst/\abs{x}$. 

The gradient estimate in \cref{eq:decay-pointwise} is now a consequence of the Schauder estimates for elliptic operators in divergence form with $\CS^{0, \alpha}$ coefficients (see \cite[Theorem 2.28]{fernandez-real_regularitytheoryellipticpde_2022}). Applying \cref{eq:boch1} to the harmonic function $u$ and plugging in \cref{eq:decay-pointwise}, \eqref{eq:decayL2} follows. 
\end{proof}

We now get a refined asymptotic $L^2$-decay estimate for the $\R^3$-Laplacian of $u$ as a consequence of \eqref{eq:decayL2}.
\begin{lemma}
Let $0 < \tau  \leq 1$ and let $(\R^3,g)$ be a $\CS^1_\tau$-Asymptotically Flat Riemannian manifold.
Then, the capacitary potential $u$ \cref{eq:p-cap-potential} of a bounded set $\Omega$ with smooth boundary satisfies
\begin{equation}
    \label{eq:refinedLapl}
    \int\limits_{A_{R, kR}} \abs{\Delta_{\eum} u}^2 \dd\mu \leq \frac{\kst}{R^{3 + 2\tau}},
\end{equation}
for $R$ large enough, where $A_{R,kR}=\set{R \leq \abs{x}\leq kR}$ and $\Delta_{\eum}$ is the Laplacian of flat $\R^3$.
\begin{proof} Let $(x^1, x^2, x^3)$ be a choice of asymptotically flat coordinates on $(M,g)$. We have
\begin{align}
    \label{eq:decompositionlapl}
    0 = \Delta_g u& = g^{ij}\partial_i \partial_j u  - \Gamma^k_{ij}\partial_ku= (g^{ij} - \eum^{ij})\partial_i \partial_j + \eum^{ij}\partial_i \partial_ju - \Gamma^k_{ij}\partial_k u\\
    &=(g^{is}\eta_{sk}g^{kj}+ O(\abs{\eta}^2))\partial_i \partial_ju + \eum^{ij}\partial_i \partial_ju - \Gamma^k_{ij}\partial_k u.
\end{align}
It follows that
\begin{equation}
    \label{eq:pointwisedecaylap}
    (\Delta_{\eum} u(x))^2 \leq \kst\left( \frac{1}{\abs{x}^{2\tau}} \abs{\D\D u (x)}^2 + \frac{1}{\abs{x}^{2 + 2\tau}} \abs{\D u}\right).  
\end{equation}
Integrating \cref{eq:pointwisedecaylap} on $A_{R, kR}$ and applying \cref{lem:decayL2}, \cref{eq:refinedLapl} follows.
\end{proof}
\end{lemma}
We are finally ready to show the pointwise asymptotic behaviour of $u$ and $\D u$ with the explicit deficit.
\begin{proposition}
\label{prop:decaywithdeficit}
Let $0 < \tau  \leq 1$ and let $(\R^3,g)$ be a $\CS^1_\tau$-Asymptotically Flat Riemannian manifold. 
Then, the capacitary potential $u$ \cref{eq:p-cap-potential} of a bounded set $\Omega$ with smooth boundary satisfies
\begin{align}
    \label{eq:decaywithdeficit}
    \abs{ u(x) - \frac{\ncapa_2(\partial \Omega)}{\abs{x}}}=O(\abs{x}^{-1-\tau}), &&  
\abs{\partial_j u +\ncapa_2(\partial \Omega) \frac{\eum_{jk}x^k}{\abs{x}^3}} =O(\abs{x}^{-2-\tau}),
\end{align}
for every $j=1,2,3$, as $\abs{x} \to +\infty$.
\end{proposition}
\begin{proof} Fix a flat chart at infinity $(M\smallsetminus K, (x^1,x^2,x^3))$ for the manifold $(M,g)$. The function $u$ can be considered as a smooth function on the set $\R^3 \smallsetminus \set{\abs{x} \leq R_0}$ which is diffeomorphic to $M \smallsetminus K$. With abuse of notation, we interpret $u$ as a smooth function on the whole $\R^3$ by extending it smoothly in the interior of $\set{\abs{x} \leq R_0}$. From now on, we are working on $\R^3$ only. Let now $v:\R^3 \to \R$ be the function
\begin{equation}
\label{eq:representation}
    v(x) = -\frac{1}{4\pi} \int_{\R^3} \frac{\Delta_{\eum} u(y)}{\abs{x- y}}  \dif y.
\end{equation}
It is well-known and direct to check that $\Delta_{\eum} v= \Delta_{\eum} u$. The main part of the proof consists in showing that  $v$ satisfies the $\CS^0$-expansion at infinity claimed for $u$ in \cref{eq:decaywithdeficit}. Hence, Liouville's theorem for harmonic functions implies that $u = v$ and that $u$ decays as stated. 

We consider any point $x$ with $\abs{x} = 2R$, and we split the integral \cref{eq:representation} into the three terms
\begin{equation}
    \label{eq:threeterms}
    \begin{split}
    v(x) = \begin{multlined}[t]-\frac{1}{4\pi}\left[\int\limits_{\set{\abs{x-y} \leq R}} \frac{\Delta_{\eum} u(y)}{\abs{x- y}} \dif y+ \int\limits_{\set{\abs{y}\wedge\abs{x-y}\geq R}}\frac{\Delta_{\eum} u(y) }{\abs{x- y}} \dif y + \int\limits_{\set{\abs{y}\leq R}}\frac{\Delta_{\eum} u(y)}{\abs{x- y}} \dif y \right] \end{multlined}
    \end{split}
\end{equation}
\smallskip
The first term is estimated through H\"older's inequality and \cref{eq:refinedLapl} as
\begin{equation}
    \label{eq:firstestimate}
    \begin{split}
    \abs{\,\int\limits_{\set{\abs{x-y} \leq R}} \frac{ \Delta_{\eum} u(y)}{\abs{x- y}} \dif y} &\leq \left(\,\smashoperator[r]{\int\limits_{\set{\abs{x-y} \leq R}} }\abs{\Delta_{\eum} u(y)}^2 \dif y\right)^{\frac{1}{2}} \left(\,\int\limits_{\set{\abs{x-y} \leq R}} \frac{1}{\abs{x - y}^2} \dif y\right)^{\frac{1}{2}} \\
    &\leq \left(\,\smashoperator[r]{\int\limits_{A_{R, 3R}}} \abs{\Delta_{\eum} u(y)}^2 \dif y\right)^{\frac{1}{2}}\left(\,\int\limits_{\set{\abs{y} \leq R}} \frac{1}{\abs{y}^2} \dif y\right)^{\frac{1}{2}}
    \leq \frac{\kst}{R^{1+\tau}},
    \end{split}
\end{equation}
where $A_{R, 3R}= \set{R \leq \abs{x} \leq 3 R}$.

To estimate the second one, since $\abs{x-y} \geq R$, we have
\begin{equation}
\label{eq:secondestimatepartial}
   \abs{\,\int\limits_{{\set{\abs{y}\wedge\abs{x-y}\geq R}}} \frac{ \Delta_{\eum} u(y) }{\abs{x- y}}\dif y} \leq \frac{\kst}{R} \int\limits_{\set{\abs{y}\geq R}} \abs{\Delta_{\eum} u(y)} \dif y.  
   \end{equation}
   Cover now $\set{\abs{y} \geq R}$ with annuli $A_j$ of the form $\set{2^jR \leq \abs{x} \leq 2^{j+1} R}$, $j\in\N$, to obtain
   \begin{equation}
       \label{eq:coverbyballs1}
       \begin{split}
       \int\limits_{\set{\abs{y} \geq R}} \abs{\Delta_{\eum} u(y)} \dif y &\leq \sum_{j = 0}^{+\infty} \int_{A_j} \abs{\Delta_{\eum} u(y)} \dif y \leq
       \sum_{j = 0}^{+\infty} \left(\int_{A_{j}} \abs{\Delta_{\eum} u(y)}^2 \dif y\right)^{\frac{1}{2}} \abs{A_{j}}^{\frac{1}{2}} \\&\leq 
       \frac{\kst}{R^{\tau}} \sum_{j = 0}^{+\infty} \frac{1}{2^{\tau j}} \leq \frac{\kst}{R^{\tau}}.
       \end{split}
   \end{equation}
   Plugging it into \cref{eq:secondestimatepartial}, it yields
   \begin{equation}
       \label{eq:secondestimate}
        \abs{\,\int\limits_{\set{\abs{y}\wedge\abs{x-y}\geq R}} \frac{\Delta_{\eum} u(y)}{\abs{x- y}}  \dif y} \leq \frac{\kst}{R^{1 +\tau}}.
   \end{equation}
   The third term in \cref{eq:threeterms} splits into
   \begin{equation}
     \label{eq:thirdtermsplit}  
     \int\limits_{\set{\abs{y}\leq R}} \frac{ \Delta_{\eum} u(y)}{\abs{x- y}} \dif y = \int\limits_{\set{\abs{y}\leq R}} \left(\frac{1}{\abs{x- y}} -\frac{1}{\abs{x}}\right) \Delta_{\eum} u(y) \dif y + \int\limits_{\set{\abs{y}\leq R}} \frac{\Delta_{\eum} u(y)}{\abs{x}}  \dif y.
   \end{equation}
The first integral on the right-hand side can be further split into
 \begin{equation}
 \label{eq:thirdtermsplit1}
 \begin{split}
    \int\limits_{\set{\abs{y} \leq R}} 
     \left(\frac{1}{\abs{x- y}} -\frac{1}{\abs{x}}\right) \Delta_{\eum} u(y) \dif y = & \int\limits_{\set{\abs{y}\leq R_0}} 
     \left(\frac{1}{\abs{x- y}} -\frac{1}{\abs{x}}\right) \Delta_{\eum} u(y) \dif y \\
     &+  \int\limits_{\set{ R_0 \leq \abs{ y } \leq R}}
     \left(\frac{1}{\abs{x- y}} -\frac{1}{\abs{x}}\right) \Delta_\delta u(y) \dif y,
     \end{split}
 \end{equation}
 for some fixed $R_0 \in (0, R)$. The first term in the right-hand side of \cref{eq:thirdtermsplit1} is easily estimated by
 \begin{equation}
     \label{eq:estimateR0}
     \abs{\,\int\limits_{\set{\abs{y} \leq R_0}} \left(\frac{1}{\abs{x- y}} -\frac{1}{\abs{x}}\right) \Delta_{\eum} u(y) \dif y} \leq \kst\int\limits_{\set{\abs{y} \leq R_0}}  \frac{\abs{y}}{\abs{x}^2} \Delta_{\eum} u(y) \leq  \frac{\kst}{R^2},
 \end{equation}
 recalling that $\abs{x} = 2R$. The second term in \cref{eq:thirdtermsplit1} is more delicate. We cover $B_{R}(0) \smallsetminus B_{R_0}(0)$ with annuli $A_j$ of the form $\set{2^jR_0 \leq \abs{x} \leq 2^{j+1} R_0}$, to obtain
 \begin{equation}
\label{eq:estimatedifficile}     
\abs{\,\int\limits_{\set{ R_0 \leq \abs{ y } \leq R}} \left(\frac{1}{\abs{x- y}} -\frac{1}{\abs{x}}\right) \Delta_{\eum} u(y) \dif y} \leq \frac{\kst}{R^2} \sum_{j=0}^{j_R} \int_{A_j}\abs{y} \abs{\Delta_\delta u(y)} \dif y,
 \end{equation}
 where $j_R$ is such that $R \in (2^{j}R_0, 2^{j+1}R_0)$. Applying H\"older's inequality and \cref{eq:refinedLapl} similarly to \cref{eq:coverbyballs1}, we get
 \begin{equation}
 \label{eq:estimatedifficile1}
 \begin{split}
     \sum_{j=0}^{j_R} \int_{A_j} \abs{y} \abs{\Delta_{\eum} u(y)} \dif y &\leq \kst \sum_{j=0}^{j_R} \left(\,\int_{A_j} \abs{\Delta_{\eum} u}^2(y)  \dif y\right)^{\frac{1}{2}} 2^{\frac{5}{2}j} R_0^{\frac{5}{2}} \\
     & \leq \kst R_0^{1-\tau}  \sum_{j=0}^{j_R} 2^{(1-\tau)j} \leq \kst\frac{ R_0^{1 - \tau}  2^{(1-\tau)j_R} 2^{1-\tau} - R_0^{1-\tau}}{2^{1-\tau} - 1}.
     \end{split}
 \end{equation}
By our choice of $j_R$ we have $R_0 2^{j(1-\tau)} \leq R^{1-\tau}$. Hence, we get
 \begin{equation}
     \label{eq:difficileversofine}
      \sum_{j=0}^{j_R} \int_{A_{j}} \abs{y} \abs{\Delta_{\eum} u(y)} \dif y \leq \kst \frac{R^{1-\tau} 2^{1-\tau} - R_0^{1-\tau}}{2^{1-\tau} - 1} = \kst R^{1-\tau}\frac{2^ {1-\tau} -\left(\frac{R_0}{R}\right)^{1-\tau}}{2^{1-\tau} - 1}. 
 \end{equation}
 It only remains to estimate the last term in \cref{eq:thirdtermsplit}. Now, $0 < \tau \leq 1$ and $R_0$ is a fixed number smaller than $R$. Plugging \cref{eq:difficileversofine} into \cref{eq:estimatedifficile}, we infer
 \begin{equation}
     \label{eq:difficilefinito}
     \abs{{\,\int\limits_{\set{R_0 \leq \abs{y} \leq R}} \left(\frac{1}{\abs{x- y}} -\frac{1}{\abs{x}}\right) \Delta_{\eum} u(y) \dif y}} \leq \frac{\kst }{R^{1+\tau}}.
 \end{equation}
 Putting together \cref{eq:threeterms}, \cref{eq:firstestimate}, \cref{eq:secondestimate}, \cref{eq:thirdtermsplit}, \cref{eq:thirdtermsplit1} and \cref{eq:difficilefinito}, we have showed that
 \begin{equation}
     \abs{{v(x) - \frac{1}{2R}\int\limits_{\set{\abs{y}\leq R}}\Delta_{\eum} u(y) \dif y}} \leq \frac{\kst}{R^{1 + \tau}},
 \end{equation}
 for any $x$ such that $\abs{x} = 2R$. By \cref{eq:coverbyballs1}, we immediately deduce that
 \begin{equation}
     \label{eq:asymptoticsenzacap}
      \abs{{v(x) -\frac{a}{\abs{x}}}} \leq \frac{\kst}{R^{1 + \tau}}
 \end{equation}
 for some $a \in \R$. In particular, $v(x) \to 0$ as $\abs{x} \to +\infty$. Since $\Delta_{\eum} v= \Delta_{\eum} u$, Liouville's theorem for harmonic functions implies that $v = u$. By the already recalled Schauder estimates \cite[Theorem 2.28]{fernandez-real_regularitytheoryellipticpde_2022}, we also have
 \begin{equation}
     \label{eq:asygradientsenzacap}
     \abs{{\partial_j u(x) +a\frac{\eum_{jl}x^l}{\abs{x}^3}}} \leq \frac{1}{R^{2 + \tau}}.
 \end{equation}
 Since $u$ is harmonic in $M \smallsetminus \overline{\Omega}$ for the metric $g$, Divergence Theorem yields
 \begin{equation}\label{eq:entercapacity}
0 = \smashoperator[r]{\int\limits_{\set{\abs{x} \leq R} \smallsetminus \overline{\Omega}}} \Delta u \dif \mu = \int\limits_{\partial \Omega} \abs{\D u} \dif \sigma + \smashoperator{\int\limits_{\set{\abs{x}=R}}} \ip{\nabla u|\nu } \dif\sigma = 4\pi \ncapa_2(\partial \Omega)+ \smashoperator{\int\limits_{\set{\abs{x} = R}}} \ip{ \nabla u| \nu } \dif \sigma.       \end{equation}
 Letting $R \to +\infty$ and comparing with \cref{eq:asygradientsenzacap}, we obtain $a=\ncapa_2(\partial\Omega)$. Plugging it into \cref{eq:asymptoticsenzacap} and \cref{eq:asygradientsenzacap} completes the proof. 
\end{proof}

Finally, a direct consequence of \cref{lem:boch} is the following integral decay estimate for the second derivatives of $u$.
\begin{proposition}
\label{cor:decayhess}
Let $0 < \tau  \leq 1$ and let $(\R^3,g)$ be a $\CS^1_\tau$-Asymptotically Flat Riemannian manifold. Then, the capacitary potential $u$ \cref{eq:p-cap-potential} of a bounded set $\Omega$ with smooth boundary satisfies
\begin{equation}
    \label{eq:decayhesswithdeficit}
\int\limits_{A_{R, kR}} \abs{ \partial_i \partial_j u +\frac{\ncapa_2(\partial \Omega)}{\abs{x}^3}\left(\eum_{ij} - 3 \frac{\eum_{ik}\eum_{jl} x^kx^l}{\abs{x}^2}\right)}^2 \dif\mu  =O(R^{-3-2\tau})
\end{equation}
for every $i,j=1,2,3$ as $R\to +\infty$, where $A_{R,kR}= \set{ R \leq \abs{x} \leq kR}$. Moreover, there exists a monotonically increasing sequence of $(t_n)_{n\in \N}$ divergent to $+\infty$ such that 
\begin{equation}\label{eq:decay_level_sets}
\int\limits_{\partial \Omega_{t_n}} \abs{\partial_i \partial_j u+\frac{\ncapa_2(\partial \Omega)}{\abs{x}^3}\left(\eum_{ij} - 3 \frac{\eum_{ik}\eum_{jl} x^kx^l}{\abs{x}^2}\right)}^2\dif \sigma =O(e^{(-4-2\tau)t_n})
\end{equation}
for every $i,j=1,2,3$ as $n \to +\infty$, where $\Omega_t= \set{ u \geq \ee^{-t}}$.
\end{proposition}
\begin{proof}
Just set $\psi = u - \ncapa_2(\partial \Omega) /\abs{x}$, and plug it into \cref{eq:boch1} in terms of $g$. Hence, we can directly exploit $\Delta u = 0$. \cref{eq:decayhesswithdeficit} directly follows using the decay estimate with deficit for the gradient obtained in \cref{prop:decaywithdeficit} and observing that 
\begin{equation}
\abs{\Delta_g \abs{x}^{-1}} = O( \abs{x}^{-3 -\tau}).
\end{equation}
By \cref{prop:decaywithdeficit}, $\set{2/R \leq u \leq 3/R}$ is contained in  $A_{R,4R}$ for some large $R$. By the Mean Value Theorem, there exists $t_n \in [R/3,R/2]$ such that
\begin{equation}
    \frac{1}{R}\int\limits_{\partial \Omega_{t_n}} \abs{\partial \partial \psi}^2\dif \sigma = \int\limits_{2/R}^{3/R} \int\limits_{\set{u=s}} \abs{\partial \partial \psi}^2 \dif \sigma\dif s \leq \kst  \int\limits_{A_{R,4R}} \abs{ \partial \partial \psi}^2\abs{\D u } \dif \sigma\leq \frac{ \kst }{R^{5+2\tau}}
\end{equation}
where we used the just proved \cref{eq:decayhesswithdeficit}. Invoking again \cref{prop:decaywithdeficit}, $R$ can be chosen large enough so that $\abs{\D u}\leq \kst_3/R^2$ which concludes since $t_n\geq R/3$.
\end{proof}

\subsection{Boundedness of the \texorpdfstring{$2$}{2}-Hawking mass}

We aim to prove that the $2$-Hawking mass of a domain with connected boundary homologous to $\partial M$ is controlled by the $\mathrm{ADM}$ mass. To ease up the computations, we start by comparing it with the classical Hawking mass on large level sets of the capacitary potential \cref{eq:p-cap-potential}. The following lemma partially inverts \cref{thm:hawking_p-hawking} and suggests that the two versions of the Hawking mass are equivalent at infinity. A similar result has been first obtained in the proof of \cite[Lemma 2.5]{agostiniani_riemannianpenroseinequalitynonlinear_2022}, under stronger asymptotic assumptions. 

\begin{lemma}\label{prop:2-hawking_vs_hawking}
    Let $(M,g)$ be a complete, noncompact $\CS^{1}_{\tau}$-Asymptotically Flat Riemannian $3$-manifold, $\tau>1/2$, with nonnegative scalar curvature and possibly with smooth and closed boundary. Assume that $H_2(M, \partial M; \Z)=\set{0}$. Let $\Omega\subseteq M$ be  bounded with connected $\CS^{1,\alpha}$-boundary homologous to $\partial M$ with $\h \in L^2(\partial \Omega)$. Then, along the sequence $(t_n)_{n\in \N}$ such that \eqref{eq:decay_level_sets} holds, we have
    \begin{equation}\label{eq:2-hawking_vs_hawking}
        \lim_{t \to +\infty} \ma^{\tph}_{H}(\partial \Omega_{t}) = \lim_{n \to +\infty} \frac{\ncapa_2(\partial \Omega_{t_n})}{32\pi} \left( 16 \pi - \int_{\partial \Omega_{t_n}}\kern-.3cm \H^2 \dif \sigma \right),
    \end{equation}
    where $\Omega_t = \set{w_2 \leq t}$.
\end{lemma}

\begin{proof} First, observe that even if $\partial \Omega$ is not smooth, we can apply \cref{prop:decaywithdeficit,cor:decayhess}. Indeed, the function $u$ is smooth away from $\partial \Omega$, by the classic elliptic regularity theory. Therefore, Sard's theorem implies that almost every level $\partial \Omega_t$ is smooth. We will work with the function $u = \ee^{-w_2}$, hence $\Omega_t= \set{u \geq \ee^{-t}}$. Let $\nu$ and $\nu_{\eum}$ be the $g$-unit normal vector field and the $\eum$-unit normal vector field to the level set $\partial \Omega_t$ respectively. Denote $\eta_{ij}=g_{ij}- \eum_{ij}$, with $\abs{\eta}=O(\abs{x}^{-\tau})$ and $\abs{\partial \eta}=O(\abs{x}^{-1-\tau})$ as $\abs{x} \to +\infty$. Moreover, denote $\psi = u - \ncapa_2/\abs{x}$, where $\ncapa_2=\ncapa_2(\partial \Omega)$. By computations
\begin{align}
    \partial_i u = -\ncapa_2 \frac{\delta_{ij} x^j}{\abs{x}^3} + \partial_i \psi,  &\qquad&  \partial_i \partial_j u =-\frac{\ncapa_2}{\abs{x}^3}\left(\eum_{ij} - 3 \frac{\eum_{ik}\eum_{jl} x^kx^l}{\abs{x}^2}\right)+   \partial_i \partial_j \psi,
\end{align}
where $\psi =O_1(\abs{x}^{-1-\tau})$ by \cref{prop:decaywithdeficit}. Using \cref{eq:inverse_metrics}, one gets
\begin{equation}\label{eq:gradient_estimate_nonsharp}
   \abs{ \partial u}_{\eum}^2 = \eum^{ij} \partial_i u\partial_j u = \abs{ \D u }_g ^2+ \eta_{ij}\eum^{ik}\eum^{jl} \partial_l u \partial_k u  + O(\abs{x}^{-4 -2 \tau}).
\end{equation}
Since the normal unit vector fields can be expressed as
\begin{align}
    \nu^j = -\frac{ g^{ij} \partial_j u}{\abs{\D u}_g}&& \text{and} && \nu_{\eum}^j = -\frac{ \eum^{ij}\partial_j u }{\abs{\partial u}_{\eum}},
\end{align}
they are related by
\begin{equation}\label{eq:unit-normal-expansion}
\begin{split}
    \nu^i &= -\frac{g^{ij} \partial_j u}{ \abs{\D u}_g}= -\left(\eum^{ij} -\eum^{ik}\eta_{lk}\eum^{jl} +O(\abs{x}^{-2\tau})\right)\frac{\abs{ \partial u}_{\eum}}{ \abs{ \D u }_g} \frac{\partial_j u}{\abs{\partial u}_{\eum}} \\ &= -\left(\eum^{ij} -\eum^{ik}\eta_{lk}\eum^{jl}+O(\abs{x}^{-2\tau})\right)\left(1+ \frac12\eta_{sm}\nu^s \nu^m + O(\abs{x}^{-2\tau})\right) \frac{\partial_j u}{\abs{\partial u}_{\eum}} \\
    &= \nu_{\eum}^i +\frac12\eta_{ks}\nu_{\eum}^k \nu_{\eum}^s\nu_{\eum}^i - \delta^{ij}\eta_{lj} \nu_{\eum}^l +O(\abs{x}^{-2\tau}).
    \end{split}
\end{equation}
Employing again \cref{prop:decaywithdeficit}, we get
\begin{align}\label{eq:hessian}
\begin{split}
    \frac{ \D_i \D_j u}{\abs{\D u}_g} &= \frac{ \partial_i \partial_j u - \Gamma^k_{ij} \partial_k u}{\abs{\partial u}_{\eum}}\frac{\abs{\partial u}_\delta}{\abs{\D u }_g}\\
    &=\begin{multlined}[t]\frac{ \partial_i \partial_j u}{ \abs{ \partial u}_{\eum}} + \frac12\eta_{kl} \nu_{\eum}^k \nu_{\eum}^l\frac{ \partial_i \partial_j u}{ \abs{ \partial u}_{\eum}} + \frac{1}{2}( -\partial_k \eta_{ij}+\partial_i \eta_{kj} +\partial_j \eta_{ki}) \nu_{\eum}^k\\+ O(\abs{x}^{-1-2\tau})+ O(\abs{x}^{2-2\tau} \abs{\partial \partial \psi}_{\delta}). \end{multlined}
    \end{split}
\end{align}
To compute the mean curvature of each level, we need to take the trace of the above quantity with respect to the metric $g^{\top}$ induced on $\partial \Omega_t$ by $g$. In particular, we need the expansion
\begin{equation}\label{eq:induced_metric}
    g_{\top}^{ij} = \eum_{\top}^{ij}- \delta_{\top}^{ik} \eta_{kl} \delta_{\top}^{lj} + O(\abs{x}^{-2\tau}),
\end{equation}
where $g_{\top}^{ij}= g^{ij}- \nu^i \nu^j $ and $\eum_{\top}^{ij}= \eum^{ij}- \nu_{\eum}^i \nu_{\eum}^j$. \cref{eq:induced_metric} follows from the expansion of $g^{ij}$ and $\nu^i$ in \cref{eq:inverse_metrics,eq:unit-normal-expansion}. Therefore, we get
\begin{equation}\label{eq:H_expansion}
\begin{split}
    \H =- g_{\top}^{ij}\frac{ \D_i \D_j u}{\abs{\D u}_g} = \begin{multlined}[t]\H_{\eum} + \frac12{\eta_{kl} \nu_{\eum}^k \nu_{\eum}^l} \H_{\eum}- \eum_{\top}^{ij}\left( \partial_j \eta_{ik}- \frac{1}{2} \partial_k \eta_{ij}\right)\nu_{\eum}^k + \eum_{\top}^{ik}\eum_{\top}^{jl} \eta_{kl}\frac{\partial_i \partial_j u}{\abs{ \partial u}_{\eum}}\\ + O(\abs{x}^{-1-2\tau})+ O(\abs{x}^{2-2\tau} \abs{\partial \partial \psi}_\eum),\end{multlined}
\end{split}
\end{equation}
where $\H$ and $\H_{\eum}$ are the mean curvatures of the level sets with respect to $g$ and $\eum$ respectively. Appealing to \cref{prop:decaywithdeficit} to estimate $\partial \psi$, we get
\begin{equation}\label{eq:u_expansion_refined}
\begin{gathered}
    \H_{\eum}= -\eum^{ij}_{\top}\frac{\partial_i \partial_j u}{\abs{\partial u}_{\eum}}= \frac{2}{\abs{x}} + O(\abs{x}^{-1-\tau})+ O(\abs{x}^{2}\abs{\partial \partial \psi}_{\delta}).
\end{gathered}
\end{equation}
Plugging \cref{eq:u_expansion_refined} into \cref{eq:H_expansion}, we finally obtain
\begin{equation}\label{eq:zzdecay_meancurvatureing}
    \H = \frac{2}{\abs{x}}+ O(\abs{x}^{-1 -\tau}) + O(\abs{x}^2 \abs{\partial \partial \psi}_\delta).
\end{equation}

On the other hand, \cref{prop:decaywithdeficit} yields
\begin{equation}\label{eq:zzdecay_grad_vs_meancurvature}
    \frac{\abs{\D u}_g}{u} = \frac{1}{\abs{x}} + O(\abs{x}^{-1 -\tau}).
\end{equation}
Therefore, \cref{eq:zzdecay_grad_vs_meancurvature},\cref{eq:zzdecay_meancurvatureing} and Young's inequality imply
\begin{equation}
    \left( \frac{\H}{2}- \frac{\abs{\D u}_g}{u} \right)^2 =O(\abs{x}^{-2-2\tau}) + O(\abs{x}^{4}\abs{\partial \partial \psi}^2_\delta).
\end{equation}
Let $(t_n)_{n\in \N}$ be the sequence given by \cref{cor:decayhess}. Using \cref{cor:decayhess} together with \cref{prop:decaywithdeficit} and \cref{lem:capacitygrowth}, we have 
\begin{equation}
    \int_{\partial \Omega_{t_n}}\left( \frac{\H}{2}- \frac{\abs{\D u}_g}{u} \right)^2 \dif \sigma = \int_{\partial \Omega_{t_n}}O(\abs{x}^{-2-2\tau})+ O(\abs{x}^{4}\abs{\partial \partial \psi}^2) \dif\sigma = O(\ee^{-2 \tau t_n}).
\end{equation}
Recalling the relation between $U_2$ and $\ma^{\tph}_H$ in \cref{eq:literature_mass_relations}, we conclude that
\begin{align}
    \lim_{t\to +\infty} \frac{8\pi}{\ncapa_2(\partial \Omega)}\ma_{H}^{\tph}(\partial \Omega_t)&=\lim_{n \to +\infty}\frac{8\pi}{\ncapa_2(\partial \Omega)}\ma_{H}^{\tph}(\partial \Omega_{t_n}) = \lim_{n \to +\infty}  U_2(\ee^{t_n})\\&= \lim_{n\to +\infty} \ee^{t_n} \left(4 \pi - \int_{\partial \Omega_{t_n}} \frac{\H^2}{4} \dif \sigma \right) + \ee^{t_n}\int_{\partial \Omega_{t_n}}\left( \frac{\H}{2}- \frac{\abs{\D u}_g}{u} \right)^2 \dif \sigma\\
    &= \lim_{n\to +\infty} \ee^{t_n} \left(4 \pi - \int_{\partial \Omega_{t_n}} \frac{\H^2}{4} \dif \sigma \right) + O(\ee^{(1- 2\tau)t_n}).
\end{align}
Since $\tau>1/2$, \cref{eq:2-hawking_vs_hawking} is proved.
\end{proof}

\medskip

The following result constitutes the conclusion of all the linear potential theoretic analysis carried out so far. 
\begin{theorem}\label{thm:bound_2_hawking}
    Let $(M,g)$ be a complete, noncompact $\CS^{1}_{\tau}$-Asymptotically Flat Riemannian $3$-manifold, $\tau>1/2$, with nonnegative scalar curvature and possibly with smooth and closed boundary. Assume that $H_2(M, \partial M; \Z)=\set{0}$. Let $\Omega\subseteq M$ be bounded with connected $\CS^{1,\alpha}$-boundary homologous to $\partial M$ with $\h \in L^2(\partial \Omega)$. Then,
    \begin{equation}\label{eq:U_2-asymptotic}
        \ma^{\tph}_H (\partial \Omega) \leq \ma_{\ADM}.
    \end{equation}
\end{theorem}
\begin{proof} For computational simplicity, we work once again with $U_2$. As in \cref{prop:2-hawking_vs_hawking}, even if $\partial \Omega$ is not smooth, we can apply \cref{cor:decayhess,prop:decaywithdeficit}. If we show that the right-hand side in \cref{eq:2-hawking_vs_hawking} is bounded by the $\mathrm{ADM}$ mass, the conclusion follows. Indeed, by the monotonicity of $ \ma^{\tph}_H$ in \cref{thm:AMMO_monotonicity} we conclude 
\begin{equation}
     \ma^{\tph}_H(\partial \Omega)\leq \lim_{t \to +\infty}  \ma^{\tph}_H(\partial \Omega_t) \leq \ma_{\ADM}.
\end{equation}
We will use the same notation of \cref{prop:2-hawking_vs_hawking} and part of the computations carried out there.

 Raising to the square \cref{eq:H_expansion} and using the Young's inequality, we get
\begin{align}\label{eq:H2-asymptotic}
\begin{split}
    \H^2 =\begin{multlined}[t] \H_{\eum}^2+ \eta_{kl} \nu_{\eum}^k \nu_{\eum}^l \H_{\eum}^2 -  \H_{\eum}\eum_{\top}^{ij}\left( 2\partial_j \eta_{ik}-\partial_k \eta_{ij}\right)\nu_{\eum}^k+ 2 \H_{\eum} \eum_{\top}^{ik}\eum_{\top}^{jl}\eta_{kl} \frac{\partial_i \partial_j u}{\abs{ \partial u}_{\eum}}\\ + O(\abs{ x}^{-2-2\tau})+ O(\abs{x}^{4-2\tau}\abs{\partial \partial \psi}_{\delta}^2) .\end{multlined}
\end{split}
\end{align}
Denoting $\dd \sigma$ and $\dd \sigma_{\eum}$ the area measures induced on the level sets by the metric $g$ and $\eum $ respectively, we have
\begin{equation}\label{eq:area-asymptotic}
    \dd \sigma = \left( 1+ \frac{1}{2}g_{\top}^{ij} \eta_{ij} + O (\abs{x}^{-2\tau})\right)\dd\sigma_{\eum}.
\end{equation}
Coupling \cref{eq:H2-asymptotic,eq:area-asymptotic}, we get the following expression of the Willmore energy of the level sets
\begin{equation}\label{eq:integrated_H2_expansion}
\begin{split}
    \int_{\partial \Omega_t}\kern-.2cm \H^2 \dif \sigma = \int_{\partial \Omega_t}\kern-.3cm\begin{multlined}[t]\H_{\eum}^2 +\eta_{kl}\nu_{\eum}^k \nu_{\eum}^l \H_{\eum}^2 +\frac{1}{2}\eum_{\top}^{kl}\eta_{kl}\H_{\eum}^2-  \H_{\eum} \eum_{\top}^{ij}\left( 2\partial_j \eta_{ik}- \partial_k \eta_{ij}\right)\nu_{\eum}^k\\+ 2 \H_{\eum} \eum_{\top}^{ik}\eum_{\top}^{jl}\eta_{kl}\frac{\partial_i \partial_j u}{\abs{ \partial u}_{\eum}} + O(\abs{x}^{-2-2\tau}) + O(\abs{x}^{4-2\tau}\abs{\partial \partial \psi}_{\delta}^2)\dif \sigma_{\eum}.\end{multlined}
    \end{split}
\end{equation}
Plugging \cref{eq:u_expansion_refined} into \cref{eq:integrated_H2_expansion} and using Young's inequality, one has
\begin{equation}\label{eq:precise_integral_H2_expansion}
\begin{split}
    \int_{\partial \Omega_t}\kern-.2cm \H^2 \dif \sigma = \int_{\partial \Omega_t}\kern-.3cm\begin{multlined}[t] \H_{\eum}^2 + \frac{4}{\abs{x}^2} \eta_{ij}\nu_{\eum}^i \nu_{\eum}^j- \frac{2}{\abs{x}^2} \eum^{kl}_{\top} \eta_{kl}- \frac{2}{\abs{x}} \eum_{\top}^{ij}\left( 2\partial_j \eta_{ik}- \partial_k \eta_{ij}\right)\nu_{\eum}^k\dif \sigma_{\eum}\\+ \int_{\partial \Omega_t}O(\abs{x}^{-2-2\tau})+ O(\abs{x}^{4}\abs{\partial \partial \psi}^2_{\delta} ) \dif \sigma_{\eum}. \end{multlined}
\end{split}
\end{equation}
Integration by parts on level sets gives
\begin{equation}
    \int_{\partial \Omega_t} \frac{1}{ \abs{x}}\H_{\eum} \nu_{\eum}^i \eta_{ik}\nu_{\eum}^k \dif \sigma_{\eum} = \int_{\partial \Omega_t} \frac{1}{ \abs{x}}\eum_{\top}^{ij} \partial_i \eta_{jk}\nu_{\eum}^k - \frac{1}{ \abs{x}}\eum_{\top}^{ij}\eta_{jk}\eum_{\top}^{kl} \frac{\partial_i \partial_l u}{\abs{\partial u}_{\eum}} \dif \sigma_{\eum},
\end{equation}
that, again by \cref{eq:u_expansion_refined} and Young's inequality, reduces to
\begin{equation}\label{eq:precise_integration_byparts}
\begin{split}
    \int_{\partial \Omega_t}\frac{2}{\abs{x}^2}\nu_{\eum}^i \eta_{ik}\nu_{\eum}^k \dif \sigma_{\eum} = \int_{\partial \Omega_t} \begin{multlined}[t][.5\textwidth]\frac{1}{\abs{x}}\eum_{\top}^{ij} \partial_i \eta_{jk}\nu_{\eum}^k + \frac{1}{\abs{x}^2}\eum_{\top}^{ij}\eta_{ks}\dif \sigma_{\eum}\\+\int_{\partial \Omega_t} O(\abs{x}^{-2-2\tau}) + O(\abs{x}^{4}\abs{\partial \partial \psi}^2 ) \dif \sigma_{\eum}.\end{multlined}
\end{split}
\end{equation}
By \cref{eq:precise_integral_H2_expansion}, \cref{eq:precise_integration_byparts}, \cref{lem:capacitygrowth} and the Willmore's inequality on $\R^n$, we conclude
\begin{align}
    \int_{\partial \Omega_t}\kern-.2cm \H^2 \dif \sigma &=\begin{multlined}[t]\int_{\partial \Omega_t}\kern-.2cm  \H_{\eum}^2 \dif \sigma_{\eum}-\int_{\partial \Omega_t}\frac{2}{\abs{x}}g^{ij}_{\top} (\partial_j g_{ik}- \partial_k g_{ij}) \nu^k \dif \sigma   \\+\int_{\partial \Omega_t}O(\abs{x}^{-2-2\tau}) + O(\abs{x}^{4}\abs{\partial \partial \psi}^2) \dif\sigma\end{multlined}\\
    &\geq\begin{multlined}[t] 16 \pi-\frac{2}{\ncapa_2(\partial \Omega_t)}\int_{\partial \Omega_t}g^{ij}_{\top} (\partial_j g_{ik}- \partial_k g_{ij}) \nu^k \dif \sigma   \\+\int_{\partial \Omega_t}O(\abs{x}^{-2-2\tau}) + O(\abs{x}^{4}\abs{\partial \partial \psi}^2) \dif\sigma.\end{multlined}
\end{align}
Along the sequence $(t_n)_{n \in \N}$ such that \cref{eq:decay_level_sets} holds, \cref{cor:decayhess} yields
\begin{align}
    \int_{\partial \Omega_{t_n}}O(\abs{x}^{-2-2\tau}) + O(\abs{x}^{4}\abs{\partial \partial \psi}^2) \dif\sigma = O(\ee^{-2\tau t_n}).
\end{align}
By \cref{lem:capacitygrowth} and $\tau>1/2$, we infer
\begin{equation}
    \lim_{n\to +\infty}\frac{\ncapa_2(\partial \Omega_{t_n})}{32 \pi }\left( 16 \pi - \int_{\partial \Omega_{t_n}}\kern-.3cm  \H^2 \dif \sigma \right) \leq\ma_{\ADM},
\end{equation}
then \cref{prop:2-hawking_vs_hawking} concludes the proof. 
\end{proof}

\subsection{From the Isoperimetric mass to the \texorpdfstring{$\ADM$}{ADM} mass}
We now have all the necessary tools to prove \cref{thm:penrose_adm_intro,thm:mass_equivalence}. For the sake of completeness, we first state the Positive Mass Theorem under optimal decay assumptions. 
\begin{theorem}
\label{thm:pmt_adm}
Let $(M,g)$ be a complete $\CS^{1}_{\tau}$-Asymptotically Flat Riemannian $3$-manifold, $\tau >1/2$, with nonnegative scalar curvature and possibly with smooth, closed and minimal boundary. Then, it holds
\begin{equation}\label{eq:adm_pmt}
0 \leq \ma_{\ADM},
\end{equation}
with the equality satisfied if and only if $(M, g)$ is isometric to flat $\R^3$.
\end{theorem}

\begin{proof}[Proof of \cref{thm:penrose_adm_intro,thm:mass_equivalence,thm:pmt_adm}]
By \cref{thm:topology} and \cref{rem:C1-foliation}, we can assume that $M$ is diffeomorphic to $\R^3$, with a finite number of balls possibly removed.
It suffices to  show that
\begin{equation}\label{eq:hawking_bound_with_adm}
    \ma_{H}(\partial \Omega) \leq  \ma_{\ADM}
\end{equation}
holds for every $\Omega \subseteq  M$ with connected smooth boundary such that any connected component of $\partial M$ is either contained in $\Omega$ or disjoint from $\Omega$. This directly implies \cref{thm:penrose_adm_intro} and \cref{thm:pmt_adm}. Moreover, by \cref{thm:bound_hawking_to_bound_iso}, we also infer $\ma_{\mathrm{iso}} \leq \ma_{\mathrm{ADM}}$. The reverse inequality is a consequence of \cite[Corollary 2.3]{fan_largespheresmallspherelimitsbrownyork_2009} (see \cref{rmk:FST_infinite_mass}), which only requires $\CS^1$-asymptotic assumptions on the metric. This result states that
\begin{equation}
    \ma_{\ADM} =\lim_{r \to +\infty} \ma_{\iso} ( \set{ \abs{x} \leq r})\leq \ma_{\iso},
\end{equation}
concluding the proof of \cref{thm:mass_equivalence}

We now turn to the proof of \cref{eq:hawking_bound_with_adm}.
By \cref{prop:jumping_to_homologous}, we can assume that $\Omega$ has $\CS^{1,\alpha}$-boundary homologous to $\partial M$ and with $\h \in L^2(\partial \Omega)$. If $\ma_{\ADM}=+\infty$, the statement would be trivially true. By \cref{thm:hawking_p-hawking,thm:bound_2_hawking}, there exists a geometric positive constant $\kst>0$ such that
\begin{equation}\label{eq:boundedness_of_mass}
    \ma_{H}(\partial \Omega) \leq \kst \ma_{H}^{\tph}(\partial \Omega)\leq \kst \ma_{\ADM}.
\end{equation}
In particular, \cref{thm:bound_hawking_to_bound_iso} implies $\ma_{\iso} \leq \kst \ma_{\ADM}<+\infty$.

Evolve now $\Omega$ by weak IMCF $w_1$ and denote $\Omega_t =\set{w_1 \leq t}$. Since the Hawking mass is monotone \cref{thm:Geroch_monotonicity_formula} and bounded by \cref{eq:boundedness_of_mass}, we have that
\begin{equation}
    \abs{\int_{\partial \Omega_t} \H^2 \dif \sigma - 16 \pi} \leq \kst \ee^{-\frac{t}{2}}
\end{equation}
for some constant $\kst$ depending only on $\partial \Omega$ and $\ma_{\ADM}$. On the other hand, integrating \cref{eq:geroch-monotonicity} we get
\begin{equation}
    \int_0^{+\infty}\ee^{-\frac{t}{2}}\int_{\partial \Omega_t}  \vert \mathring{\h}\vert^2 \dif \sigma \dif t <+\infty.
\end{equation}
In conclusion, there exists a sequence of $t_n$'s, diverging at $+\infty$ as $n\to +\infty$, such that 
\begin{align}
\label{eq:sequence}
 \int_{\partial \Omega_{t_n}}\vert\mathring{\h}\vert^2 \dif \sigma \leq  \kst  \ee^{-\frac{t_n}{2}} && \text{ and }&&\int_{\partial \Omega_{t_n}} \vert\h\vert^2 \dif \sigma \leq  \kst  .
\end{align}
We now want to prove that a similar behaviour also holds for the quantities expressed with respect to the flat metric $\eum$. Again we denote $\eta_{ij} = g_{ij}- \eum_{ij}$. Let $\nu$ and $\nu_\eum$ be the unit normal vector fields to $\partial \Omega_t$ with respect to the metric $g$ and $\eum$ respectively. Let $\omega_i = g_{ij} \nu^j$ and $\omega^\eum_i = \eum_{ij} \nu_\eum^j$ be their respective duals. As in the proof of \cref{prop:2-hawking_vs_hawking} (see also \cite[(7.7)]{huisken_inversemeancurvatureflow_2001}), we have 
\begin{align}\label{eq:vectors_convergence}
    \omega^\eum_i = \omega_i + O(\abs{x}^{-\tau}), && \nu^i_\eum = \nu^i + O(\abs{x}^{-\tau}), && \abs{\omega_\eum}_g = 1  + O(\abs{x}^{-\tau}).
\end{align}
Denote $\h^\eum$ the second fundamental form of $\partial \Omega_t$ with respect to the metric $\eum$. By \cite[(7.10)]{huisken_inversemeancurvatureflow_2001} and \cref{eq:vectors_convergence}, we know that
\begin{equation}
    \h_{ij} - \h^\eum_{ij} =-\Gamma^k_{ij} \omega_k^\eum + (1- \vert\omega_\eum\vert_g) \h_{ij} = O(\abs{x}^{-1- \tau}) + O(\abs{\h} \abs{x}^{-\tau}),
\end{equation}
Here, $\Gamma_{ij}^k$ denotes the Christoffel symbols of the metric $g$. Taking the trace of it with respect to $\eum^{ij}_\top= \eum^{ij} - \nu_\eum^i \nu^j_\eum$, we also get that
\begin{equation}
    \H- \H_\eum = - \omega^\eum_k \Gamma^k_{ij}\eum^{ij}_\top + (1- \vert{\omega_\eum}\vert_g) \h_{ij}\eum^{ij}_\top + (g^{ij}_\top - \eum^{ij}_{\top}) \h_{ij}=O(\abs{x}^{-1- \tau}) + O(\abs{\h} \abs{x}^{-\tau}),
\end{equation}
where $g_{\top}^{ij} = g^{ij}- \nu^i\nu^j$ and $\H_\eum$ is the mean curvature of $\partial \Omega_t$ with respect to the metric $\eum$. Recalling that $w_1(x)= O(2\log\abs{x})$ by \cref{eq:boundimcf}, \cref{eq:area-asymptotic} and \cref{eq:sequence}, we get  
\begin{equation}
    \int_{\partial \Omega_{t_n}} \vert \mathring{\h}_\eum\vert_{\eum}^2 \dif \sigma_\eum \leq \kst \int_{\partial \Omega_{t_n}} \vert \mathring{\h}\vert^2 + O(\abs{x}^{-2-2\tau}) + O(\abs{\h}^2 \abs{x}^{-2\tau}) \dif \sigma \leq \kst \left( \ee^{-\frac{t}{2}} + \ee^{-\tau t}\right)
\end{equation}
We are then in position to apply \cite[Theorem 1.1]{delellis_nearlyumbilical_2005} and deduce that
\begin{equation}
    \label{eq:delellis}
    \int_{\partial \Omega_{t_n}} \abs{{ \H_\eum -  \sqrt{\frac{16\pi}{\abs{ \partial \Omega_{t_n}}_\eum}}}}^2 \dif \sigma_\eum \leq \kst \ee^{-\frac{t_n}{2}}.
\end{equation}
Arguing as above and employing again \cref{eq:area-asymptotic}, this estimate holds also when the quantities are referred to the metric $g$. Therefore, we have 
\begin{align}\label{eq:sharp_mean_curvature_HI}
     \int_{\partial \Omega_{t_n}} \abs{{ \H -  \sqrt{\frac{16\pi}{\abs{ \partial \Omega_{t_n}}}}}}^2 \dif \sigma \leq \kst \ee^{-\frac{t_n}{2}}.
 \end{align}
Moreover, by \cref{eq:sequence} and \cref{eq:sharp_mean_curvature_HI} it also holds
\begin{align}\label{eq:sharp_second_ff_HI}
    \int_{\partial \Omega_{t_n}} \abs{{\h-\sqrt{\frac{4\pi}{\abs{ \partial \Omega_{t}}}}  g^{\top}}}^2 \dif \sigma \leq\int_{\partial \Omega_{t_n}}\vert\mathring{\h}\vert + \frac{\sqrt{2}}{2}\abs{{ \H -  \sqrt{\frac{16\pi}{\abs{ \partial \Omega_{t}}}}}} \dif \sigma \leq  \kst \ee^{-\frac{t_n}{2}}.
\end{align}
Following the same computations carried out for \cref{prop:2-hawking_vs_hawking} (see also \cite[Asymptotic Comparison Theorem 7.4]{huisken_inversemeancurvatureflow_2001}) and employing \cref{eq:boundimcf,eq:sharp_mean_curvature_HI,eq:sharp_second_ff_HI}, one obtains
\begin{equation}\label{eq:main_inequality_HI}
\begin{split}
    16 \pi \ma_{H}(\partial \Omega_{t_n})&=\begin{multlined}[t]\sqrt{ \frac{\abs{\partial \Omega_{t_n}}}{16\pi}} \int_{\partial \Omega_{t_n}}\kern-.1cm 2 \H g_{\top}^{im}\eta_{ml}g_{\top}^{lj} \h_{ij}-\frac{1}{2} \H^2 g_{\top}^{ij}\eta_{ij}- \H^2 \nu^i\nu^j\eta_{ij}\dif \sigma \\ +\sqrt{ \frac{\abs{\partial \Omega_{t_n}}}{16\pi}} \int_{\partial \Omega_{t_n}}2 \H g_{\top}^{ij} \nu^l \left(\D_i \eta_{jl}-\D_l \eta_{ij}\right) +O(\abs{\h}^2 \abs{x}^{-2\tau})
    O(\abs{x}^{-2-2\tau})\dif \sigma\end{multlined}\\
    & \leq\begin{multlined}[t] \int_{\partial \Omega_{t_n}}\kern-.1cm g^{ij}(\partial_i g_{jl} - \partial_l g_{ij} ) \nu^l \dif\sigma+ \kst \ee^{\frac{1}{4}(1-2\tau)t_n}+ \kst\ee^{\frac{1}{2}(1- 2 \tau) t_n}.
    \end{multlined}
    \end{split}
\end{equation}
Since $\tau>1/2$, we finally have
\begin{equation}
    \lim_{t \to +\infty} \ma_{H}(\partial \Omega_{t}) = \lim_{n\to +\infty} \ma_{H}(\partial \Omega_{t_n}) \leq \lim_{n\to +\infty} \frac{1}{16 \pi} \int_{\partial \Omega_{t_n}}\kern-.1cm g^{ij} (\partial_i g^{\top}_{jl} - \partial_l g^{\top}_{ij}) \nu^l \dif \sigma = \ma_{\ADM}.\qedhere
\end{equation}
\end{proof}

\begin{remark}\label{rmk:FST_infinite_mass}
    In \cite{fan_largespheresmallspherelimitsbrownyork_2009}, the authors assume that the scalar curvature belongs to $L^1(M)$. Here, we do not assume any \emph{a priori} integrability for $\mathrm{R}_g$. Nonetheless, $\ma_{\ADM}$ is still a well-defined geometric invariant since the scalar curvature is nonnegative. However, it could be infinite, in which case the computations in \cite{fan_largespheresmallspherelimitsbrownyork_2009} are meaningless. One can adapt the arguments to show that $\ma_{\ADM} = +\infty$ implies $\ma_{\iso} = +\infty$. Assume $\ma_{\ADM}=+\infty$ and take $\ma\in \R$. All computations in \cite[Lemma 2.2]{fan_largespheresmallspherelimitsbrownyork_2009} still work, until passing to the limit the quantity in the definition of the $\ADM$ mass. Since $\ma < \ma_{\ADM}=+\infty$, one can replace \cite[(2.12)]{fan_largespheresmallspherelimitsbrownyork_2009} with
    \begin{equation}
        \frac{\dd}{\dd r }\abs{\set{\abs{x}=r}} \leq\frac1r{\abs{\set{\abs{x}=r}}}+4 \pi r+ \int_{\set{\abs{x}=r}}(g_{ij}-\delta_{ij})\frac{x^ix^j}{r^3} \dif \sigma_\delta - 8\pi \ma
    \end{equation}
    for $r$ large enough. Following the computations in \cite[Theorem 2.2]{fan_largespheresmallspherelimitsbrownyork_2009} one obtains
    \begin{equation}
        \abs{\set{\abs{x} \leq r}}\geq\frac{1}{2} r \abs{\set{\abs{x}=r}} - \frac{2\pi}{3}r^3 + 2 \pi \ma r^2.
    \end{equation}
    Plugging it into the definition of the Isoperimetric mass, we conclude
    \begin{equation}
        \ma_{\iso} \geq \limsup_{r \to +\infty} \ma_{\iso}(\set{\abs{x} \leq r}) \geq \ma
    \end{equation}
    for every $\ma\in \R$.
\end{remark}

\begin{appendices}
\crefalias{section}{appendix}
\section{Generalised de l'H\^opital's rule}
For the sake of completeness, we report here the de l'H\^opital's rule we crucially used in this paper. The proof is inspired to \cite[Theorem II]{taylor_hospitalrule_1952}.

\begin{theorem}[de l'H\^opital's rule]\label{thm:delhopital}
Let $I\subseteq \R$ be an open and possibly unbounded interval of the real line and $a\in \R \cup\set{\pm \infty}$ one of its endpoints. Let $f,g:I \to \R$ be two functions satisfying the following properties:
\begin{enumerate}
    \item $f$ is locally absolutely continuous on $I$;
    \item $g$ is everywhere differentiable on $I$ and $g(x)\neq 0$ and $g'(x) \neq 0$ for every $x \in I$.
    \item  there exists a sequence $(t_k)_{k \in \N}$, $t_k \in I$, $t_k\to a\in I$ as $k\to +\infty$ such that one of this two conditions is satisfied
    \begin{enumerate}
        \item\label{item:condition_1}$f(t_k)\to 0$ and $g(t_k)\to 0$ when $k \to +\infty$, or
        \item\label{item:condition_2} $\abs{g(t_k)} \to +\infty$ as $k \to +\infty$.
    \end{enumerate}
\end{enumerate}
Then,
\begin{equation}\label{eq:generalised_del'hopital}
    \liminf_{x\to a} \frac{f'(x)}{g'(x)}\leq \liminf_{k\to +\infty}\frac{f(t_k)}{g(t_k)} \leq \limsup_{k\to +\infty} \frac{f(t_k)}{g(t_k)} \leq \limsup_{x \to a} \frac{f'(x)}{g'(x)}.
\end{equation}
\end{theorem}

\begin{proof} Without loss of generality we can assume that $a$ is the right endpoint and thus $t_k \leq a$ for every $k \in \N$. By Darboux's theorem, $g'$ is of constant sign, then $g(x) \neq g(y)$ for every $x,y \in I$. Since $f$ is absolutely continuous we have that
\begin{equation}\label{eq:del'hopital_equivalence}
    \frac{f(y)-f(x)}{g(y)-g(x)}=\frac{ \int_x^y f'(t) \dif t}{g(y)-g(x)} = \frac{\int_x^y g'(t) \frac{f'(t)}{g'(t)} \dif t }{g(y)-g(x)} 
\end{equation}
for every $x,y \in I$. Let now $y \in (x,a)$. By the fundamental theorem of calculus applied to the function $g$, we directly get that
\begin{equation}\label{eq:del'hopital_estimates}
     \inf_{s \in (x,a)} \frac{f'(s)}{g'(s)} \leq\frac{\int_x^y g'(t) \frac{f'(t)}{g'(t)} \dif t}{g(y)-g(x)} \leq \sup_{s \in (x,a)} \frac{f'(s)}{g'(s)}
\end{equation}
hold for every $x<y<a$ and $x,y\in I$. Coupling together \cref{eq:del'hopital_equivalence,eq:del'hopital_estimates}, we obtain
\begin{align}\label{eq:del'hopital_main1}
    \inf_{s \in (x,a)} \frac{f'(s)}{g'(s)} &\leq\frac{\dfrac{f(y)}{g(x)}-\dfrac{f(x)}{g(x)}}{\dfrac{g(y)}{g(x)}-1} \leq \sup_{s \in (x,a)} \frac{f'(s)}{g'(s)}
    \intertext{and} \label{eq:del'hopital_main2}
    \inf_{s \in (x,a)} \frac{f'(s)}{g'(s)} &\leq\frac{\dfrac{f(x)}{g(y)}-\dfrac{f(y)}{g(y)}}{1-\dfrac{g(x)}{g(y)}} \leq \sup_{s \in (x,a)} \frac{f'(s)}{g'(s)}.
\end{align}
If \cref{item:condition_1} holds, choosing $y=t_k$ and sending $k\to +\infty$ in \cref{eq:del'hopital_main1} we conclude
\begin{equation}
    \inf_{s \in(x,a)} \frac{f'(s)}{g'(s)} \leq\frac{f(x)}{g(x)} \leq \sup_{s \in (x,a)} \frac{f'(s)}{g'(s)}.
\end{equation}
\cref{eq:generalised_del'hopital} follows from it, choosing again $x=t_k$ and sending $k \to +\infty$. The case \cref{item:condition_2}  follows similarly, using \cref{eq:del'hopital_main2} accordingly. \end{proof}

\section{Stable two-sided bounds for the \texorpdfstring{$p$}{p}-capacitary potential}
The study of the existence of a (weak) proper IMCF in increasingly broader contexts has been a challenging topic since its introduction by Huisken and Ilmanen \cite{huisken_inversemeancurvatureflow_2001}. In their setting, i.e. in $\CS^1$-Asymptotically Flat case, the existence follows by building an explicit subsolution of the form $u_\alpha= (2-\alpha) \log \abs{x}$. Indeed, at large distances, this function satisfies
\begin{equation}
    \div \left( \frac{\D u_\alpha}{\abs{\D u_\alpha}}\right)> \abs{\D u_\alpha}.
\end{equation}
Hence, the level sets of $u_\alpha$ are strictly mean convex and move faster than the IMCF would dictate. A consequence of this property is that IMCF exists and is controlled from below by $u_\alpha$ up to an additive constant. Similarly, one can produce an upper bound by proving that $(2+\alpha)\log\abs{x}$ is a supersolution.
The major issue with this procedure is that it appears to heavily rely on the $\CS^1$ properties of the metric, even if the functional minimised by the weak IMCF does not. 

Conversely, the recent approach by Xu \cite{xu_isoperimetrypropernessweakinverse_2023} shows that a proper solution to weak IMCF exists in a wide framework, including $\CS^0$-Asymptotically Flat manifolds. This approach provides a lower bound akin to that obtained through the subsolution $u_\alpha$ in the $\CS^1$-Asymptotically Flat setting. However, it is not clear to us how to provide a similar upper bound, that we actually needed for the proof of \cref{thm:asysvarchild}.

We follow here the strategy put forth by Moser \cite{moser_inversemeancurvatureflow_2007}, that we used in \cref{prop:c01nonpar}. Indeed, the IMCF can be approximated in the limit as $p\to 1^+$ by the solutions $w_p$ of 
\begin{equation}\label{eq:pIMCF}
 \begin{cases}
 \Delta_p w_p &=& \abs{ \D w_p}^p & \text{on $M\smallsetminus \overline{\Omega}$,}\\
 w_p &=& 0 & \text{on $\overline{ \Omega}$,}\\
 w_p(x) &\to& +\infty & \text{as $\abs{x} \to +\infty$},
 \end{cases}
\end{equation}
where $\Delta_p f= \div(\abs{\D f}^{p-2} \D f)$ is the usual $p$-Laplacian. The function $u_p = \ee^{-w_p/(p-1)}$ is $p$-harmonic. As a key tool in the proof, we use a Harnack inequality for positive $p$-harmonic functions on coordinate spheres $\set{\abs{x} \leq R}$ with a constant that remains stable as $p\to 1^+$ and is independent of $R$. To the best of the authors' knowledge, a proof of the Harnack inequality for the IMCF without passing through the nonlinear potential theory is unavailable in the literature. 

\begin{lemma}\label{lem:p-harnack}
    Let $(M, g)$ be a complete $\CS^0$-Asymptotically Flat Riemannian $3$-manifold and $1<p<3$. Then, there exists a constant $\kst_H>0$ independent of $p$ such that for every $R>0$ large enough one has
    \begin{equation}\label{eq:p-harnack}
         \sup_{\set{\abs{x}=R}} u_p(x) \leq \kst_H^{\frac{1}{p-1}} \inf_{\set{\abs{x} =R}} u_p(x),
    \end{equation}
    for every nonnegative function $u_p$ which is $p$-harmonic in $\set{R/2 \leq \abs{x} \leq 3R/2}$.
\end{lemma}
\begin{proof}
    We denote $B(z,r)= \set{\abs{x-z} < r}$ for some $z \in \set{\abs{x}=R}$. Given a $p$-harmonic functions $u_p$ defined in $\set{R/2 \leq \abs{x} \leq 3R/2}$ of a given Riemannian manifold $(M, g)$, it satisfy the Harnack inequality 
    \begin{equation}
    \label{eq:harnack}
        u_p(x) \leq \kst_H^{\frac{1}{p-1}}(B(z, R/2)) u_p(y)
    \end{equation}
    for any $x, y \in B(z, R/4)$ for a constant $\kst_H$ that depends only on the Poincar\'e, Sobolev and Alfhors constants of $B(z, R/2)$. This result is obtained in \cite{rigoli_noteonpsubharmonicfunctionscomplete_1997}, see \cite[Theorem 3.4 and Remark 3.5]{mari_flowlaplaceapproximationnew_2022} for a discussion of the explicit constants. By asymptotic flatness, for $R$ large enough we can cover $\set{\abs{x} = R}$ with a uniform number $N$ of coordinate balls $B(z, R/6)$. As $\kst_H$ remains controlled in the $\CS^0$-asymptotic flat regime, we deduce \cref{eq:p-harnack}.
\end{proof}

Assuming that a solution to \cref{eq:pIMCF} exists, $u_p= \ee^{-w_p/(p-1)}$ is $p$-harmonic outside $\Omega$. Hence, \cref{eq:p-harnack} applies and shows that
\begin{equation}
    \inf_{\set{\abs{x}=R}} w_p(x) \geq \sup_{\set{\abs{x}=R}}w_p(x) - \kst_H.
\end{equation}
Therefore, establishing stable lower and upper bounds for $w_p$ at any point of $\set{\abs{x}=R}$ is sufficient to control the whole behaviour of the function. Such bounds follow exploiting the properties of $p$-capacities and the maximum principle. Observe that the lower bound is also the keystone to establish the existence of a solution to \eqref{eq:pIMCF}.

\begin{proposition}\label{thm:stronglyp-li-yau}
    Let $(M, g)$ be a complete $\CS^0$-Asymptotically Flat Riemannian $3$-manifold. Then, for any $\Omega\subseteq M$ with $\CS^{1,\alpha}$-boundary homologous to $\partial M$ there exists a solution to the problem \cref{eq:pIMCF} for every $1<p<3$. Moreover,
    \begin{equation}\label{eq:li-yau-p}
        (3-p) \log\abs{x} - \kst_1 \leq w_p(x) \leq (3-p) \log \abs{x} +\kst_2
    \end{equation}
    for constants $\kst_1, \kst_2>0$  not depending on $p$.
\end{proposition}
\begin{proof}
    We first show the existence of a solution to \cref{eq:pIMCF}. To achieve it, we first approximate \cref{eq:pIMCF} with $w_p^R = -(p-1)\log u_p^R$, where $u_p^R$ solves
\begin{equation}\label{eq:p-cap-potential-inproof}
    \begin{cases}
     \Delta_p u_p^R& =& 0 &\text{on $\set{\abs{x} < 2R} \smallsetminus \overline{\Omega}$,}\\
     u_p^R&=&1 & \text{on $\partial \Omega$,}\\
     u_p^R&=& 0 & \text{on $\set{\abs{x}= 2R}$.}
    \end{cases}
\end{equation}
For a compact subset $K \subset A$ in an open set $A$, we define the (normalized) relative $p$-capacity of $K$ in $A$ as
\begin{equation}\label{eq:p-capacity-inproof}
    \ncapa_p(K,A) = \inf \set{\frac{1}{4\pi}\left(\frac{p-1}{3-p}\right)^{p-1}\int_{A} \abs{\D v}^p \dif \mu\st v \in \CS_c^\infty(A),\, v \geq 1\text{ on } K}.
\end{equation}
It is well known (see e.g. \cite[Lemma 3.8]{holopainen_nonlinearpotentialtheoryquasiregular_1990}) that 
\begin{equation}
\label{eq:capexp}
    \ncapa_p(\set{w_p^R\leq t}, \set{\abs{x} < 2R }) = \ee^{t} \ncapa_p(\partial \Omega,  \set{\abs{x} < 2R }).
\end{equation}
Let now $M_R = \max\set{w_p^R(x)\st \abs{x} =R}$. By the maximum principle, we have
$\set{\abs{x} \leq R} \subset \set{w_p^R \leq M_R}$. By the monotonicity of $p$-capacities with respect to inclusion and \cref{eq:capexp}, we thus get
\begin{equation}
    \ncapa_p(\set{\abs{x} = R}, \set{\abs{x} < 2R }) \leq \ncapa_p(\set{w_R^p \leq M_r}, \set{\abs{x} < 2R }) = \ee^{M_R} \ncapa_p(\partial \Omega,  \set{\abs{x} < 2R }),
\end{equation}
yielding
\begin{equation}
\label{eq:stimasup}
   M_R \geq  \log\left(\frac{\ncapa_p(\set{\abs{x} = R}, \set{\abs{x} < 2R })}{\ncapa_p(\partial \Omega)}\right).
\end{equation}
As a consequence of \cite[Theorem 3.2]{pigola_connectivityinfinitymanifoldq_2014}, $\CS^0$-Asymptotically Flat manifolds support a global isoperimetic inequality. Thus, we can employ it to infer the $p$-isocapacitary inequality (see e.g. \cite[(7)]{grigoryan_isoperimetricinequalities_1998})
\begin{equation}
\label{eq:piso}
  \ncapa_p(\set{\abs{x} = R}, \set{\abs{x} < 2R }) \geq \kst_p \abs{\set{\abs{x} \leq  R}}^{\frac{3-p}{3}}.
\end{equation}
The constant $\kst_p$ tends to the isoperimetric constant as $p\to 1^+$, and in particular can be chosen independently of $p$. Moreover, $\ncapa_p(\partial \Omega)$ tends to $(4\pi)^{-1}\abs{\partial \Omega^*}$ as $p \to 1^+$ by \cite[Theorem 1.2]{fogagnolo_minimisinghullscapacityisoperimetric_2022}. Plugging \cref{eq:piso} into \cref{eq:stimasup}, and exploiting the asymptotic flatness to estimate the volume of big geodesic balls we are thus left with
\begin{equation}
\label{eq:stimasupexpl}
    M_R \geq (3-p)\log(R) - \kst
\end{equation}
By the Harnack inequality \cref{lem:p-harnack}, recalling that $w_p^R = -(p-1) \log u_p^R$, we have
\begin{equation}
\label{eq:lowerbound_w_pR}
w_p^R(x) \geq (3-p) \log \abs{x} - \kst.
\end{equation}
On the other hand, by \cite[Theorem 1.1]{kotschwar_localgradientestimatesharmonic_2009} the functions $w_p^R$ satisfy a gradient bound on each compact set of $\set{\abs{x} < 2R} \smallsetminus \overline{\Omega}$.  As in the proof of \cite[Proposition 3.3]{kotschwar_localgradientestimatesharmonic_2009}, one can construct a smooth barrier function $v$ for any $x_0 \in \partial \Omega$. This function is independent of $p$, satisfies $v(x_0) = 0$, and $w_p^R \leq v$ in a neighbourhood of $x_0$. Coupling this information with the interior gradient bound, we can pass to the limit as $R \to + \infty$ in any compact set of $M \smallsetminus \Omega$. The limit $w_p$ is a (weak) solution to $\Delta_p w_p = \abs{\D w_p}^p$ in $M\smallsetminus \overline{\Omega}$, which agrees with $0$ continuously on $\partial \Omega$. The lower bound \cref{eq:lowerbound_w_pR} is preserved in the limit as $R\to +\infty$. Hence $w_p\to + \infty$ as $\abs{x}\to +\infty$. Therefore, $w_p$ is a solution to \cref{eq:pIMCF} and in particular it satisfies
\begin{equation}
    w_p(x) \geq (3-p) \log \abs{x} -\kst_1,
\end{equation}
where $\kst_1>0$ not depending on $p$.

We are left to show the upper bound claimed in \cref{eq:li-yau-p}. Let $m_R = \min\set{w_p(x)\st \abs{x} =R}$, where $w_p$ is again the solution to \cref{eq:pIMCF}. The comparison principle guarantees that $\set{w_p \leq m_R} \subset \set{\abs{x} \leq R}$. Again by the monotonicity of $p$-capacities with respect to inclusion and \cref{eq:capexp}, we deduce that
\begin{equation}
\label{eq:chaincapacities}
\ee^{m_R} \ncapa_p(\partial \Omega) = \ncapa_p(\partial \set{w_p \leq m_R}) \leq \ncapa_p(\set{\abs{x} = R}). 
\end{equation}
Applying the Harnack inequality as above, we deduce that
\begin{equation}
\label{eq:upperboundw_p}
  \max_{y \in \set{\abs{x} =R}}w_p(y) \leq \log (\ncapa_p(\set{\abs{x} = R})) + \kst ,
\end{equation}
  for a constant $\kst$ that does not depend on $p$.  Now, the asymptotic behaviour of $p$-capacities at infinity discussed in in \cite[Lemma 2.21]{benatti_asymptoticbehaviourcapacitarypotentials_2022} implies
  \begin{equation}
      \label{eq:asy1cap}
      \lim_{R\to + \infty} \frac{\ncapa_p(\set{\abs{x} = R})}{R^{3-p}} =1,
  \end{equation}
  and that this limit is stable in $p$. Namely, the proof of \cite[Lemma 2.21]{benatti_asymptoticbehaviourcapacitarypotentials_2022} in fact shows that, given $\varepsilon > 0$, there exists $R_\varepsilon > 0$ independent of $p > 1$ such that 
  \begin{equation}
1 - \varepsilon \leq \frac{\ncapa_p(\set{\abs{x} = R})}{R^{3-p}} \leq 1 + \varepsilon
  \end{equation}
  for any $R \geq R_\varepsilon$.
  Plugging this piece of information into \cref{eq:upperboundw_p}, we infer the upper bound in \cref{eq:li-yau-p}.
\end{proof}

Consider a $\CS^0$-Asymptotically Flat complete Riemannian manifold $(M,g)$ and a point $o \in M$. Slightly modifying the above proof, one can show that a weak solution $w_p$ to 
 \begin{equation}\label{eq:pIMCFfromapoint}
\begin{cases}
 \Delta_p w_p &=& \abs{ \D w_p}^p & \text{on $M\smallsetminus \set{o}$,}\\
 w_p(x) &\to& -\infty & \text{as $\dist(x,o) \to 0$,}\\
 w_p(x) &\to& +\infty & \text{as $\dist(x,o) \to +\infty$}
 \end{cases}
\end{equation}
exists and satisfies the bounds in \cref{eq:li-yau-p} for constants $\kst_1, \kst_2$ that do not depend on $p$. Namely, it suffices to replace $u_p^R$ solving \cref{eq:p-capacity-inproof} with the $p$-Green's function $G_p^R$ of $\set{\abs{x} <  2R}$ with pole $o \in \set{\abs{x} < 2R}$. The proof runs substantially the same. The only modification consists of using a standard uniform lower bound on $G_p$ instead of the barrier argument at $\partial \Omega$ (see e.g. \cite[Corollary 2.8]{mari_flowlaplaceapproximationnew_2022}). This bound follows from the Laplace comparison and the known asymptotics of $G_p$ at the pole.

\begin{remark}
\label{rem:C01pnonpar}
Consider $(w_p)_{p >1}$ the family of solutions either of \cref{eq:pIMCFfromapoint} if $\partial M= \varnothing$, or \cref{eq:pIMCF} with $\Omega= \partial M$.
Using a gradient bound from \cite{kotschwar_localgradientestimatesharmonic_2009,mari_flowlaplaceapproximationnew_2022}, we get a local uniform limit of $(w_p)_{p >1}$ which is a weak IMCF. The weak IMCF satisfies the double bound \cref{eq:boundimcf}, hence it solves \cref{eq:IMCFfromapoint} or \cref{eq:IMCF} starting at $\partial M$. In particular, $\CS^0$-Asymptotically flat manifolds are examples of strongly $1$-nonparabolic manifolds.

\smallskip

These results never take advantage of dimension $3$, and can be carried out with obvious modifications in any dimension $n \geq 2$.
\end{remark}
\end{appendices}
\begingroup
\setlength{\emergencystretch}{1em}
\printbibliography
\endgroup

\end{document}